\documentclass[a4,twoside,10pt]{amsart}

\usepackage{a4wide}
\usepackage{amsfonts}
\usepackage{amsmath}
\usepackage{amssymb}
\usepackage{amsthm}
\usepackage[english]{babel}
\usepackage{bbm}
\usepackage{enumerate}
\usepackage{fontenc}
\usepackage[colorlinks=true,linkcolor=blue,citecolor=red,urlcolor=blue]{hyperref}
\usepackage[latin1]{inputenc}
\usepackage[scr=rsfs]{mathalpha}
\usepackage{orcidlink}
\usepackage{stackrel}
\usepackage[english]{varioref}
\usepackage[all]{xy}
\usepackage{url}

\title{Lefschetz principle-type theorems for curve semistable Higgs sheaves and applications to elliptic surfaces}

\author{Armando Capasso}
\address{Universit\`a degli Studi di Trieste, P.le Europa 1, Trieste (Italy), C.A.P. 34127}
\email{armando.capasso@units.it}

\thanks{A.C. is member of INdAM-GNSAGA. ORCID: 0009-0001-5463-7221 \orcidlink{0009-0001-5463-7221}}


\date{\today}

\theoremstyle{plain}
\newtheorem{theorem}{Theorem}[section]

\newtheorem{corollary}[theorem]{Corollary}
\newtheorem{lemma}[theorem]{Lemma}

\newtheorem{proposition}[theorem]{Proposition}

\theoremstyle{definition}
\newtheorem{conjecture}[theorem]{Conjecture}
\newtheorem{definition}[theorem]{Definition}
\newtheorem{definitions}[theorem]{Definitions}
\newtheorem{ex}[theorem]{Example}
\newenvironment{example}{\begin{ex}}{\hfill{$\triangle$}\end{ex}}
\newtheorem{rem}[theorem]{Remark}
\newenvironment{remark}{\begin{rem}}{\hfill{$\Diamond$}\end{rem}}
\newtheorem{remarkl}[theorem]{Remark}
\newtheorem{remarks}[theorem]{Remarks}
\newtheorem*{Proof}{Proof}
\newenvironment{prf}{\begin{Proof}}{\hfill\text{(Q.e.d.)}\end{Proof}}

\DeclareMathOperator{\End}{\mathrm{End}}
\DeclareMathOperator{\Gr}{\mathrm{Gr}}
\DeclareMathOperator{\hGr}{\mathfrak{Gr}}
\DeclareMathOperator{\cHom}{\mathcal{H}om}
\DeclareMathOperator{\Id}{Id}
\let\Im\relax
\DeclareMathOperator{\Im}{Im}
\DeclareMathOperator{\Pic}{Pic}
\DeclareMathOperator{\pr}{pr}
\DeclareMathOperator{\proj}{\mathrm{Proj}}
\DeclareMathOperator{\qr}{qr}
\DeclareMathOperator{\bQuot}{\mathbf{Quot}}
\DeclareMathOperator{\fQuot}{\mathfrak{Quot}}
\DeclareMathOperator{\Quot}{\mathrm{Quot}}
\DeclareMathOperator{\rad}{\mathrm{rad}}
\DeclareMathOperator{\rank}{\mathrm{rank}}
\DeclareMathOperator{\second}{\prime\prime}
\DeclareMathOperator{\spec}{\mathrm{Spec}}
\DeclareMathOperator{\specm}{\mathrm{SpecMax}}
\DeclareMathOperator{\Supp}{\mathrm{Supp}}
\DeclareMathOperator{\Tr}{Tr}
\DeclareMathOperator{\uzero}{\underline{0}}

\newcommand{\C}{\mathbb{C}}
\newcommand{\F}{\mathbb{F}}
\newcommand{\I}{\mathbb{I}}
\newcommand{\bL}{\mathbb{L}}
\newcommand{\K}{\mathbb{K}}
\newcommand{\N}{\mathbb{N}}
\newcommand{\bP}{\mathbb{P}}
\newcommand{\Q}{\mathbb{Q}}
\newcommand{\Z}{\mathbb{Z}}

\newcommand{\bfP}{\mathbf{P}}
\newcommand{\bfQ}{\mathbf{Q}}

\newcommand{\cC}{\mathcal{C}}
\newcommand{\cE}{\mathcal{E}}
\newcommand{\cF}{\mathcal{F}}

\newcommand{\cK}{\mathcal{K}}
\newcommand{\cL}{\mathcal{L}}
\newcommand{\cO}{\mathcal{O}}
\newcommand{\cQ}{\mathcal{Q}}
\newcommand{\cS}{\mathcal{S}}
\newcommand{\cU}{\mathcal{U}}
\newcommand{\cV}{\mathcal{V}}
\newcommand{\cW}{\mathcal{W}}

\newcommand{\fE}{\mathfrak{E}}
\newcommand{\wfE}{\widetilde{\mathfrak{E}}}
\newcommand{\fF}{\mathfrak{F}}
\newcommand{\fm}{\mathfrak{m}}
\newcommand{\fQ}{\mathfrak{Q}}
\newcommand{\fU}{\mathfrak{U}}

\newcommand{\rA}{\mathrm{A}}
\newcommand{\rh}{\mathrm{h}}
\newcommand{\rH}{\mathrm{H}}
\newcommand{\rJ}{\mathrm{J}}
\newcommand{\rM}{\mathrm{M}}
\newcommand{\rQ}{\mathrm{Q}}

\newcommand{\sfAlgSp}{\mathsf{AlgSp}}
\newcommand{\sfC}{\mathsf{C}}
\newcommand{\osfM}{\overline{\mathsf{M}}}
\newcommand{\sfQuot}{\mathsf{Quot}}
\newcommand{\sfSch}{\mathsf{Sch}}
\newcommand{\sfSet}{\mathsf{Set}}

\newcommand{\ocC}{\overline{\cC}}

\newcommand{\whC}{\widehat{C}}
\newcommand{\whcO}{\widehat{\mathcal{O}}}
\newcommand{\whL}{\widehat{L}}
\newcommand{\whS}{\widehat{S}}

\newcommand{\wtC}{\widetilde{C}}
\newcommand{\wtcC}{\widetilde{\mathcal{C}}}

\newcommand{\wtg}{\widetilde{g}}
\newcommand{\wtL}{\widetilde{L}}
\newcommand{\wts}{\widetilde{s}}
\newcommand{\wtS}{\widetilde{S}}
\newcommand{\wtgamma}{\widetilde{\gamma}}
\newcommand{\wteta}{\widetilde{\eta}}
\newcommand{\wtsigma}{\widetilde{\sigma}}
\newcommand{\wtvarphi}{\widetilde{\varphi}}

\begin{document}

\begin{abstract}
I prove ``Lefschetz principle''-type theorems for slope semistable and curve semistable Higgs sheaves on smooth projective varieties, defined over an algebraically closed field of characteristic $0$. These theorems are applied to reduce a conjecture, about curve semistable Higgs bundles, from the previous general setting to the complex case. This conjecture is equivalent to triviality of Chern classes of H-nflat Higgs bundles, which are particular curve semistable Higgs bundles. Where the base variety is a (complex) Jacobian elliptic surfaces, we prove that H-nflat Higgs bundles determine a well-defined subset of moduli spaces of Higgs sheaves: the ``H-nflat locus''. We prove that this is set-wise fixed by a natural action, and at ``limits'' one finds nilpotent, H-nflat Higgs bundles. As consequence, the conjecture for Jacobian elliptic surfaces is reduced to consider nilpotent, H-nflat Higgs bundles; moreover, their Higgs fields vanish, hence the conjecture holds (also) for Jacobian elliptic surfaces.
\end{abstract}

\maketitle

\tableofcontents

\markboth{Lefschetz principle-type theorems for curve semistable Higgs sheaves and applications to elliptic surfaces}{Armando Capasso}

\section*{Introduction}

\noindent Let $X$ be a scheme over an algebraically closed field $\K$ of characteristic $0$, unless otherwise indicated. ``Lefschetz principle''-type theorems state that some properties of $X$ (\emph{e.g.} to be a projective variety) are stable under opportune base change from $\K$ to complex numbers field $\C$ and \emph{vice versa}.\medskip

\noindent Assuming that $X$ is a smooth, projective variety, we prove a ``Lefschetz principle''-type theorem for semistable Higgs sheaves on $X$ (Lemma \ref{lem2.2} and Corollary \ref{cor2.1}) explicitly, \emph{i.e.} we prove that \emph{semistability} and \emph{curve semistability} (see Definition \ref{def0.1}) conditions for Higgs sheaves on $X$ are stable under any base field change of $\K$.\medskip

\noindent To be clear, first ``Lefschetz principle''-type theorem concerns slope (semi)stable Higgs sheaves on $X$ (Lemma \ref{lem2.2}). It is proved using \emph{Higgs-Quot schemes} (Definition \ref{def1.2}): they are \emph{moduli spaces} classifying \emph{families of quotient Higgs sheaves} of a given Higgs sheaf (Definition \ref{def1.3}), they carry on a \emph{universal family of quotient Higgs sheaves} (Theorem \ref{th1.1}) and they are well behaved under base change (Proposition \ref{prop1.1}). The best of my knowledge, these mathematical objects has been never studied in literature.\medskip

\noindent Second ``Lefschetz principle''-type theorem concerns curve semi)stable Higgs sheaves on $X$ (Corollary \ref{cor2.1}). Since there we consider the pull-backs of a Higgs sheaf on irreducible, smooth projective curves mapped to a base change of a projective variety, it is natural to consider these objects as Higgs sheaves on families of stable maps. To formalize this point of view, one uses the \emph{Higgs-Quot stack} of a fixed Higgs sheaf over the \emph{stack of stable maps into} $X$ (Definition \ref{def1.1} and Theorem \ref{th1.2}). In general, families of stable maps are \emph{algebraic spaces} over schemes; to avoid this generality, we prove that it is ever possible to cover these base schemes with an \emph{\'etale covering} so that the algebraic space obtained via base change is a projective scheme (Proposition \ref{propA.1} and \ref{propA.2}).\medskip 

\noindent We are interested in studying rank $r\geq2$ curve semistable Higgs bundles $\fE=(E,\varphi)$ over $X$, because they are semistable with respect to all polarizations of $X$. And it is conjectured by Bruzzo and Gra\~na Otero that their \emph{discriminant class} $\displaystyle\frac{1}{2r}c_2\left(E\otimes E^{\vee}\right)$ vanishes (Conjecture \ref{conj}).
\begin{remarkl}\label{rem0.2}
Assuming $\K=\C$, the best of my knowledge, Conjecture \ref{conj} has been proved in the following cases:
\begin{enumerate}[a)]
\item\label{rem0.2.a} $\varphi=0$, by  \cite[Proposition 1.3]{C:P} and \cite[Corollary 3.2]{B:B:G};
\item\label{rem0.2.b} $r=2$, by \cite[Theorems 4.5, 4.8 and 4.9]{B:GO:HR};
\item $X$ has nef tangent bundle, by \cite[Corollary 3.15]{B:LG};
\item $\dim(X)=2$ and the Kodaira dimension $\kappa(X)$ of $X$ is either $-\infty$ or $0$ (\cite[Theorem 3.6, Propositions 3.11, 3.12, Corollary 3.8]{B:LG} and \cite[Theorem 6.4]{B:L:LG}).
\item $\dim(X)=2$, $\kappa(X)=1$ and other technical hypotheses, see \cite[Proposition 5.6]{B:P};
\item $X$ is a simply-connected Calabi-Yau variety, by \cite[Theorem 4.1]{B:C};
\item if $X$ satisfies the conjecture and $Y$ is a fibred projective variety over $X$ with rationally connected fibres, then $Y$ does the same, by \cite[Proposition 3.11]{B:LG};
\item if $X$ satisfies the conjecture then any finite \'etale quotient $Y$ of $X$ does the same, by \cite[Proposition 3.12]{B:LG};
\item $\fE$ has a Jordan-H\"older filtration whose quotients are curve semistable, their first Chern class vanishes and have rank at most $2$, by \cite[Theorem 3.2, Corollary 3.3]{B:C} and previous item \ref{rem0.2.b};
\item particular Higgs bundles described in \cite{B:GO:2}. \hfill{$\Diamond$}
\end{enumerate}
\end{remarkl}
\noindent Using the ``Lefschetz principle''-type theorem for curve semistable Higgs sheaves, we reduce Conjecture \ref{conj} to complex case, \emph{i.e.} if the Conjecture holds for smooth, complex, projective varieties then it holds for any smooth projective variety defined over $\K$ (cfr. Remark \ref{rem2.1} and Theorem \ref{th2.2}).\medskip

\noindent Form another point of view, among curve semistable Higgs bundles there are some interesting one: those that, after pull-backed to curves, have first Chern classe numerically equivalent to zero (cfr. Proposition \ref{prop3.1}); these are called \emph{numerically flat Higgs bundles} (\emph{H-nflat}, for short; see Definition \ref{def3.1}). H-nflat Higgs bundles satisfy nice properties (see for example \cite{B:GO:2,B:B:G,B:C}). Moreover, Conjecture \ref{conj} is known to be equivalent to triviality of Chern classes of H-nflat Higgs bundles (\cite[Corollary 3.2]{B:B:G}); and by previous papers of Bruzzo, Lanza and Lo Giudice, in order to prove this conjecture, it is enough to study H-nflat Higgs bundles over minimal, smooth, projective surfaces (\cite[Proposition 3.13]{B:LG} and \cite[Subsection 3.1]{L:LG}).\medskip

\noindent Thus, considering the previous list, it is natural to want to study H-nflat Higgs bundles over \emph{elliptic surfaces}, and the second main result of this paper is the triviality of Chern classes of H-nflat Higgs bundles over \emph{Jacobian} elliptic surfaces (Theorem \ref{th4.1}). The idea is to reduce further Conjecture \ref{conj} to only \emph{nilpotent}, H-nflat Higgs bundles (Lemma \ref{lem3.1} and Remark \ref{remC.1}). Hence, via classical Grothendieck's and Atiyah's papers (\cite{G:A,A:MF}) on vector bundles over projective line and elliptic curves, respectively, we prove that the Higgs field of nilpotent, H-nflat Higgs bundles vanishes (Lemma \ref{lem4.3}). This allows us to give a proof of Conjecture \ref{conj} for Jacobian elliptic surfaces, \emph{i.e.} minimal, smooth, complex, projective surfaces whose projection onto the base curve has a section (Definition \ref{def4.1}).\medskip

\noindent Aside, restricting the focus to \emph{non-isotrivial} elliptic surfaces, we prove some generalisation of results proved by Bruzzo and Peragine (\cite{B:P}) about Higgs bundles over these surfaces, Franco, Garc\'ia-Prada and Newstead (\cite{F:GP:N}) about Higgs bundles over elliptic curves. Furthermore, we give an explicit example which proves the necessity of the non-isotriviality assumption (Example \ref{ex4.1}).\medskip

\noindent At the end, this paper contains a first appendix, where we collect know facts about \emph{moduli stack of stable maps to} $X$. And there we prove some basic results on families of stable maps, for which we are not able to give a reference.\medskip

\noindent In the second appendix we collect facts known about \emph{projective varieties and their base change}. There we give explicit proofs of these statements.\medskip

\noindent In the third and last appendix, we collect facts known about \emph{moduli spaces of Higgs sheaves over smooth, complex, projective surfaces}. Above all them: the natural action of $\C^{\times}$ over these spaces. This is the main mathematical tool used to reduce Conjecture \ref{conj} for Jacobian elliptic surfaces to only nilpotent, H-nflat Higgs bundles; this is also the cornerstone for the second main result of this paper.\medskip

\noindent This third appendix is necessary because H-nflat Higgs bundles are slope semistable, but it is unknown whether they are semistable (see Definition \ref{def0.2} and Lemma \ref{lem3.2}). Hence it needs to study the natural action of $\C^{\times}$ on the set of \emph{JH-equivalence classes} of H-nflat Higgs bundles over Jacobian elliptic surface (the \emph{H-nflat locus}, see Equation \eqref{eq3.1}): at first sight this set is not well-defined, but using the \emph{Principle of Induction} on the rank of Higgs sheaves, in parallel to the proof of the conjecture, one has that this set is well-defined. Furthermore, one proves that it is set-wise fixed by this action (Proposition \ref{prop3.2}); thus one proves that limits of H-nflat Higgs bundles are nilpotent in addition to be H-nflat (Corollary \ref{cor3.1}).

\subsection*{Slope (semi)stability and (semi)stability conditions for Higgs sheaves}

\noindent Let $X$ be a $\K$-scheme, let $\Omega^1_X$ be the cotangent sheaf of $X$ and let $\displaystyle\Omega^p_X=\bigwedge^p\Omega^1_X$ for $p\in\N_{\geq1}$.
\begin{definition}
A \emph{Higgs sheaf} $\fE$ on $X$ is a pair $(\cE,\varphi)$ where $\cE$ is an $\cO_X$-coherent sheaf equipped with a morphism $\varphi\colon\cE\to\cE\otimes\Omega^1_X$ called \emph{Higgs field} such that the composition
\begin{displaymath}
\varphi\wedge\varphi\colon\cE\xrightarrow{\varphi}\cE\otimes\Omega^1_X\xrightarrow{\varphi\otimes\Id}\cE\otimes\Omega^1_X\otimes\Omega^1_X\xrightarrow{\Id\otimes\pi}\cE\otimes\Omega^2_X
\end{displaymath}
vanishes. A \emph{Higgs subsheaf} of $\fE$ is a $\varphi$-\emph{invariant subsheaf} $\cF$ of $\cE$, \emph{i.e.} $\varphi(\cF)\subseteq\cF\otimes\Omega_X^1$. A \emph{quotient Higgs sheaf} of $\fE$ is a quotient sheaf of $\cE$ such that the corresponding kernel is $\varphi$-invariant. A \emph{Higgs bundle} is a Higgs sheaf on $X$ whose underlying coherent sheaf is locally free.
\end{definition}
\noindent Henceforward, let assume $X$ be projective, let $H$ be a polarization of $X$ and let $\fE=(\cE,\varphi)$ be a torsion-free Higgs sheaf on $X$, if not otherwise indicated. One defines the \emph{slope} of $\fE$ as $\displaystyle{\mu(\fE)=\frac{1}{\rank(\cE)}\int_Xc_1(\det(\cE))\cap H^{n-1}\in\Q}$; for simplicity, $\_\cap H^{n-1}$ means $\_\cap c_1\left(\cO_X(H)\right)^{n-1}$.
\begin{definition}\label{def0.1}
$\fE$ is \emph{slope} $H$-(\emph{semi})\emph{stable} if $\mu(\fF)\stackrel[(\leq)]{\textstyle<}{}\mu(\fE)$ for every Higgs subsheaf $\fF$ of $\fE$ with ${0<\rank(\fF)<\rank(\fE)}$. Or $\fE$ is \emph{slope} $H$-(\emph{semi})\emph{stable} if and only if $\mu(\fE)\stackrel[(\leq)]{\textstyle<}{}\mu(\fQ)$ for every torsion-free quotient Higgs sheaf $\fQ$ of $\fE$ with $0<\rank(\fQ)<\rank(\fE)$, equivalently. $\fE$ is \emph{slope} $H$-\emph{polystable} if it is a direct sum of slope $H$-stable Higgs sheaves having the same slope. $\fE$ is \emph{curve semistable} if $f^{\ast}\fE$ is slope semistable for any $f\colon C\to X$, where $C$ is an irreducible, smooth, projective curve. In the other eventuality, $\fE$ is $H$-\emph{unstable}.
\end{definition}
\noindent For simplicity, we shall skip any reference to the fixed polarization $H$ of $X$, if there is no confusion.\medskip

\noindent Let $\cL$ be an ample invertible sheaf on $X$; let $\cE$ and $\cF$ be coherent sheaves on $X$. It is known that for $m\gg1$ the Euler characteristic $\displaystyle\chi(\cE,\cL,m)=\sum_{i\in\N_{\geq0}}(-1)^iH^i\left(X,\cE\otimes\cL^m\right)$ of $\cE$ with respect to $\cL$ is a polynomial function (the so-called \emph{Hilbert polynomial} $p^{\cL}_{\cE}$ \emph{of} $\cE$ \emph{with respect to} $\cL$) with coefficient in $\Q$ and degree $\dim\Supp(\cE)$ (cfr. \cite[Theorem 18.6.1]{FOAG}). One posits $p^{\cL}_{\cE}\preceq p^{\cL}_{\cF}$ if $p^{\cL}_{\cE}(m)\leq p^{\cL}_{\cF}(m)$ for $m\gg1$.
\begin{definition}[{\cite[Proposition 3.3.(ii)]{HC:MG}}]\label{def0.2}
A torsion-free Higgs sheaf $\fE=(\cE,\varphi)$ over $X$ is (\emph{semi})\emph{stable with respect to} $\cL$ if the inequality
\begin{displaymath}
p^{\cL}_{\cE,norm}=\frac{p^{\cL}_{\cE}}{\rank(\cE)}\stackrel[(\preceq)]{}{\prec}\frac{p^{\cL}_{\cQ}}{\rank(\cQ)}=p^{\cL}_{\cQ,norm}
\end{displaymath}
holds for any torsion-free quotient Higgs sheaf $\fQ=\left(\cQ,\wtvarphi\right)$ of $\fE$ with $0<\rank(\cQ)<\rank(\cE)$. In the other eventuality, $\fE$ is \emph{unstable with respect to} $\cL$.
\end{definition}
\noindent For simplicity, we shall skip any reference to the fixed ample invertible sheaf $\cL$ on $X$, if there is no confusion.
\begin{remark}[{see \cite[Propositions 3.1.(i) and .(ii)]{HC:MG}}]\label{rem0.1}
The following implications are well known facts:
\begin{displaymath}
\xymatrix{
\text{slope stable} \ar@{=>}[r] & \text{stable} \ar@{=>}[r] & \text{semistable} \ar@{=>}[r] & \text{slope semistable}
}.
\end{displaymath}
The equivalences between (semi)stability and slope (semi)stability hold if $\dim(X)=1$ or ${\rank(\cE)=1}$.~\end{remark}

\noindent \textbf{Notations and conventions.} $\K$ is an algebraically closed field of characteristic $0$, unless otherwise indicated. By a projective variety $X$ we mean a projective, integral scheme over $\spec(\K)$ ($\K$-scheme, for short) of dimension $n\geq1$ and of finite type; in particular, up to isomorphisms, $X$ is a closed subscheme of $\bP^M_{\K}$ for some $M\in\N_{\geq1}$. If $n\in\{1,2\}$ we shall write projective curve or projective surface, respectively. Whenever we consider a morphism $f\colon C\to X$, we understand that $C$ is a smooth, projective curve; and $\cO_X$ understands the structure sheaf of the scheme $X$. We shall denote by $E$ a rank $r\geq1$ vector bundle over $X$, \emph{i.e.} a rank $r$ locally free sheaf of $\cO_X$-modules on $X$, while we shall use the script character $\cE$ to indicate any coherent sheaf. Somewhere we shall confuse interchangeably vector bundles and locally free sheaves.

\section{Higgs-Quot scheme and Higgs-Quot stack}\label{HiggsQuot}

\subsection{Higgs-Quot scheme}

\noindent Let $X$ be a projective $\K$-scheme and let $\fE=(\cE,\varphi)$ be a Higgs sheaf on $X$.
\begin{definition}\label{def1.3}
By \emph{family of quotient Higgs sheaves of $\fE$ over $X$ parameterised by a $\K$-scheme $S$} we mean a pair $\left(\fF=(\cF,\Psi),q\right)$ where
\begin{itemize}
\item $\fF$ is a Higgs sheaf on $X_S=X\times_{\K}S$ such that $\cF$ is flat over $S$;
\item $q$ is an epimorphism of Higgs sheaves on $X_S$ such that $q\colon\pr_1^{\ast}\fE\to\fF$, where $\pr_1\colon X_S\to X$ is the canonical projection.
\end{itemize}
Two such families $\left(\fF_1,q_1\right)$ and $\left(\fF_2,q_2\right)$ are \emph{equivalent} if $\ker\left(q_1\right)=\ker\left(q_2\right)$.
\end{definition}
\noindent It is easy to check that one can pull-back families, since flatness condition is preserved by base change; hence one has a set-valued contravariant functor $\bQuot_{\fE/X/\K}\colon\left(\sfSch_{/\K}\right)^{opp}\to\sfSet$ which maps any $\K$-scheme $S$ to $\displaystyle\left\{\begin{tabular}{l}
families of quotient Higgs sheaves\\ of $\fE$ on $X$ parameterised by $S$
\end{tabular}\right\}/\text{equivalence}$.\medskip

\noindent Let $\cL$ be a fixed ample invertible sheaf over $X$. Let $(\cF,q)$ be a family of quotient sheaves of $\cE$ on $X$ parameterised by a $\K$-scheme $S$; it is known that for any $s\in S$ the Euler characteristic $\chi\left(\cF_{\vert X_s},\cL_{\vert X_s}\right)$ is locally constant (cfr. \cite[Theorem 24.7.1]{FOAG}), hence the function ${s\in S\mapsto p_{\cF_{\vert X_s},\cL_{\vert X_s}}\!\!\in\Q[\lambda]}$ is locally constant; here $X_s=\pr_2^{-1}(s)$ and $\pr_2\colon X_S\to S$ is the canonical projection. This implies that the functor $\bQuot_{\cE/X/\K}$ naturally decomposes as a co-product:
\begin{displaymath}
\bQuot_{\cE/X/\K}=\coprod_{p\in\Q[\lambda]}\bQuot_{\cE/X/\K}^{p,\cL}
\end{displaymath}
where the functor $\bQuot_{\cE/X/\K}^{p,\cL}\colon\left(\sfSch_{/\K}\right)^{opp}\to\sfSet$ sends any $\K$-scheme $S$ to\newline
$\displaystyle\left\{\begin{tabular}{l}
families of quotient Higgs sheaves of $\fE$ on $X$ parameterised by $S$ such that\\the Hilbert polynomial of $\cF_{\vert X_s}$ with respect to $\cL_{\vert X_s}$ is $p\in\Q[\lambda]$ for any $s\in S$.
\end{tabular}\right\}/\text{equivalence}$. Each one of these functors $\bQuot_{\cE/X/\K}^{p,\cL}$ is representable by a projective $\K$-scheme $\Quot_{\cE/X/\K}^{p,\cL}$ (see for example \cite[Theorem 2.1.3, Remark 2.1.5.(5) and (6)]{A:J}).
\begin{remark}
Since each $\Quot_{\cE/X/\K}^{p,\cL}$ is a fine moduli space, there exists a \emph{universal family} $\cU_{\cE/X/\K}^{p,\cL}$ \emph{of quotient sheaves of} $\cE$ \emph{on} $X$, \emph{with fixed Hilbert polynomial} $p$ \emph{with respect to} $\cL$.
\end{remark}
\noindent Let $\cU$ be the gluing sheaf of $\cU_{\cE/X/\K}^{p,\cL}$'s, let $\pr_{\rQ}\colon X\times_{\K}\Quot_{\cE/X/\K}\to\Quot_{\cE/X/\K}$ and ${\pr_X\colon X\times_{\K}\Quot_{\cE/X/\K}\to X}$ be the canonical projections; one has the short exact sequence
\begin{displaymath}
\xymatrix{
0\ar[r] & \cK\ar[r]^(.375){\mu} & \pr_X^{\ast}\cE\ar[r]^{\epsilon} & \pr_{\rQ}^{\ast}\cU\ar[r] & 0
}
\end{displaymath}
where $\cK=\ker(\epsilon)$. Assuming also that $X$ is smooth, let $\Omega^1_X$ be the cotangent bundle of $X$; tensorising the previous short exact sequence with $\pr_X^{\ast}\Omega^1_X$, one has the following exact sequence
\begin{displaymath}
\xymatrix{
0\ar[r] & \cK\otimes\pr_X^{\ast}\Omega^1_X\ar[r]^(.45){\mu\otimes\Id} & \pr_X^{\ast}\cE\otimes\pr_X^{\ast}\Omega^1_X\ar[r]^(.475){\epsilon\otimes\Id} & \pr_{\rQ}^{\ast}\cU\otimes\pr_X^{\ast}\Omega^1_X\ar[r] & 0
}.
\end{displaymath}
One considers the morphism $\left(\epsilon\otimes\Id_{\pr_X^{\ast}\Omega^1_X}\right)\circ\psi\circ\mu\colon\cK\to\pr_{\rQ}^{\ast}\cU\otimes\pr_X^{\ast}\Omega^1_X$ and the corresponding zero locus $Z$, where $\psi\colon\pi^{\ast}_X\cE\to\pi^{\ast}_X\cE\otimes\pi^{\ast}_X\Omega^1_X$ is the pull-back of $\varphi$ via $\pr_X$.
\begin{lemma}
$\pr_{\rQ}^{\ast}\cU$ restricted to $Z$ is a quotient Higgs sheaf of $\pr_X^{\ast}\cE$.
\end{lemma}
\noindent This follows by noting that $\cK$ is a Higgs subsheaf of $\pr_X^{\ast}\cE$, if both are restricted to $Z$.\medskip

\noindent Since $\Quot_{\cE/X/\K}$ is a projective $\K$-scheme, it is proper and quasi-compact (\cite[Theorem II.4.9]{H:RC:2}) hence $\pr_{\rQ}$ is a closed and quasi-compact morphism (\cite[Proposition 10.3.(1)]{G:W:1}), $\pr_{\rQ}(Z)=\fQuot_{\fE/X/\K}$ is a closed subscheme of $\Quot_{\cE/X/\K}$ and it coincides with its scheme-theoretic image (\cite[Corollary 9.4.5]{FOAG}).
\begin{theorem}\label{th1.1}
$\fQuot_{\fE/X/\K}$ represents a closed subfunctor of $\bQuot_{\cE/X/\K}$, and ${\fU=\cU_{\vert\fQuot_{\fE/X/\K}}}$ is the \emph{universal family of quotient Higgs sheaves of} $\fE$ \emph{on} $X$, \emph{i.e.} for any $\K$-scheme $S$, for any family of quotient Higgs sheaves $(\fF,q)$ of $\fE$ parameterised by $S$, there exists a unique morphism $f\colon S\to\fQuot_{\fE/X/\K}$ such that $\left(\Id_X\times f\right)^{\ast}\left(\pr_X^{\ast}\fE\xrightarrow{\epsilon}\pr_{\fQ}^{\ast}\fU\right)=\pr_1^{\ast}\fE\xrightarrow{q}\fF$.
\end{theorem}
\begin{prf}
Let $(\fF,q)$ be a family of quotient Higgs sheaves of $\fE$ parameterised by $S$. By the \emph{Universal Property} of Quot schemes, there exists a unique morphism $f\colon S\to\Quot_{\cE/X/\K}$ such that $\left(\Id_X\times f\right)^{\ast}\left(\pr_X^{\ast}\cE\xrightarrow{\epsilon}\pr_{\cQ}^{\ast}\cU\right)=\pr_1^{\ast}\cE\xrightarrow{q}\cF$. Since $\pr_1^{\ast}\fE=\left(\Id_X\times f\right)^{\ast}\left(\pr_X^{\ast}\fE\right)$, $\ker(q)$ is a Higgs subsheaf of $\pr_1^{\ast}\fE$ and $\ker(q)=\ker\left(\left(\Id_X\times f\right)^{\ast}(\epsilon)\right)=\left(\Id_X\times f\right)^{\ast}(\ker(\epsilon))=\left(\Id_X\times f\right)^{\ast}(\cK)$, by all this it follows that
\begin{displaymath}
\ker(q)=\left(\Id_X\times f\right)^{\ast}(\cK)\xrightarrow{\left(\Id_X\times f\right)^{\ast}\left(\pr_X^{\ast}\varphi\right)}\left(\Id_X\times f\right)^{\ast}\left(\cK\otimes\pr_X^{\ast}\Omega^1_X\right)=\ker(q)\otimes\pr_1^{\ast}\Omega^1_X.
\end{displaymath}
Equivalently
\begin{displaymath}
\cK\xrightarrow{\pr_X^{\ast}\varphi}\cK\otimes\pr_X^{\ast}\Omega^1_X
\end{displaymath}
\emph{i.e.} $\cK$ is a Higgs subsheaf of $\pr_X^{\ast}\fE$; by definition $\pr_{\rQ}\fU$ is a quotient Higgs sheaf of $\pr_X^{\ast}\fE$. All this is is possible if and only if $f$ takes values in $\fQuot_{\fE/X/\K}$.
\end{prf}
\begin{definition}\label{def1.2}
$\fQuot_{\fE/X/\K}$ is called \emph{Higgs-Quot scheme} of $\fE$ over $X$.
\end{definition}
\noindent For the aims of this paper, we prove only the following property of Higgs-Quot schemes.
\begin{proposition}\label{prop1.1}
Let $S$ be a $\K$-scheme, and let $\fE_S$ be the pull-back of $\fE=(\cE,\varphi)$ over $X_S$. There exists a canonical isomorphism of schemes $S\times_{\K}\fQuot_{\fE/X/\K}\cong\fQuot_{\fE_S/X_S/S}$\footnote{
Since there exist $a,b\in\Z$ such that $\cE$ is a quotient of $\cO_X(-b)^{\oplus a}$, by base change $\cE_S$ is a quotient of $\cO_{X_S}(-b)^{\oplus a}$ and \cite[Theorem 2.1.3]{A:J} works (cfr. also \cite[Remarks 2.1.5.(5) and .(6)]{A:J}). Thus one repeats the previous constructions and finds the scheme $\fQuot_{\fE_S/X_S/S}$.
}.\end{proposition}\begin{prf}
Let $T$ be a $S$-scheme, $T$ is also a $\K$-scheme of course. Given a family $(\fF,q)$ of quotient Higgs sheaves of $\fE$ on $X$ parameterised by $T$ viewed as a $\K$-scheme, there exists a unique morphism ${f\colon T\to\fQuot_{\fE/X/\K}}$ associated to this family, as explained in Theorem \ref{th1.1}. By \cite[Exercise 7.2.2]{B:S} $f$ corresponds to a unique morphism $g_1\colon T\to S\times_{\K}\fQuot_{\fE/X/\K}$ of $S$-schemes; this determines a unique family $\left(\fF_S,q_S\right)$ of Higgs quotient sheaves of $\fE_S$ parameterised by $T$ given by $\left(\Id_{X_S}\times g_1\right)^{\ast}\left(\qr_{X_S}^{\ast}\fE_S\to\qr^{\ast}\fU\right)$, where $\qr\colon X_S\times_S\left(S\times_{\K}\fQuot_{\fE/X/\K}\right)\to\fQuot_{\fE/X/\K}$ and $\qr_{X_S}\colon X_S\times_S\left(S\times_{\K}\fQuot_{\fE/X/\K}\right)\to X_S$ are the canonical projections. Again, by Theorem \ref{th1.1}, $\left(\fF_S,q_S\right)$ corresponds to a unique morphism $g_2\colon T\to\fQuot_{\fE_S/X_S/S}$ of $S$-schemes. \emph{Vice versa}, formally $\fF_S$ is a quotient Higgs sheaf on $X_S\times_ST$, this last one is canonically isomorphic to $X\times_{\K}T$ as $\K$-scheme via a morphism $j$. Thus $j^{\ast}\left(\pr_{X_S,T}^{\ast}\fE_S\xrightarrow{q_S}\fF_S\right)$ corresponds to a unique morphism $T\to\fQuot_{\fE/X/\K}$, where $\pr_{X_S,T}X_S\times_ST\to X_S$ is a canonical projection. Again by \cite[Exercise 7.2.2]{B:S} this morphism corresponds to a unique morphism $T\to S\times_{\K}\fQuot_{\fE/X/\K}$ of $S$-schemes which determines $\left(\fF_S,q_S\right)$, \emph{i.e.} this morphism is $g_1$.\smallskip

\noindent By all this, one has a bijection $k(T)\colon\hom_S\left(T,S\times_{\K}\fQuot_{\fE/X/\K}\right)\to\hom_S\left(T,\fQuot_{\fE_S/X_S/S}\right)$ for any $S$-scheme $T$. Moreover, repeating the previous reasoning, for any $S$-schemes $T$ and $V$ and for any $\ell\in\hom_S(T,V)$, one has the following commutative diagram
\begin{displaymath}
\xymatrix{
\hom_S\left(V,S\times_{\K}\fQuot_{\fE/X/\K}\right)\ar[d]_{\hom_S\left(\ell,S\times_{\K}\fQuot_{\fE/X/\K}\right)}\ar[r]^(.525){k(V)} & \hom_S\left(V,\fQuot_{\fE_S/X_S/S}\right)\ar[d]^{\hom_S\left(\ell,\fQuot_{\fE_S/X_S/S}\right)}\\
\hom_S\left(T,S\times_{\K}\fQuot_{\fE/X/\K}\right)\ar[r]_(.525){k(T)} & \hom_S\left(T,\fQuot_{\fE_S/X_S/S}\right)
};
\end{displaymath}
by definition $\hom_S\left(\_,S\times_{\K}\fQuot_{\fE/X/\K}\right)$ and $\hom_S\left(\_,\fQuot_{\fE_S/X_S/S}\right)$ are naturally isomorphic functors, hence by \emph{Yoneda's Lemma} one has the claim.
\end{prf}
\begin{remark}\label{rem1.1}
Instead to work in the category of $\K$-schemes, one may repeat all these reasoning in the category of schemes and one finds the functors $\bQuot_{\fE/X/\Z}^{p,\cL}$, the schemes $\Quot_{\fE/X/\Z}^{p,\cL}$ and $\fQuot_{\fE/X/\Z}^{p,\cL}$ and so on.
\end{remark}

\subsection{Higgs-Quot stack over the moduli stack of stable maps}

\noindent Here we construct the \emph{Higgs-Quot stack} $\sfQuot_{\fE/\osfM_{g,m}(X,\beta)}$ \emph{over the stack of stable maps into} $X$ \emph{with fixed class} ${\beta\in\rA^{n-1}(X)}$; here $X$ is a projective $\K$-scheme and $\fE=(\cE,\varphi)$ be a Higgs sheaf on $X$. The preliminary notions about the stack $\osfM_{g,m}(X,\beta)$ are explained in Appendix \ref{ModStackStMaps}.
\begin{definition}\label{def1.1}
The objects of $\sfQuot_{\fE/\osfM_{g,m}(X,\beta)}$ are tuples $\left(\cC\to S,\wtS\to S,\cC\to X,\sigma_{\ast},\fF,q\right)$ where:
\begin{enumerate}[a)]
\item $\left(\cC\to S,\cC\to X,\sigma_{\ast}\right)$ is an object of $\osfM_{g,m}(X,\beta)$;
\item\label{def1.1.b} $\wtS\to S$ is a morphism of schemes for which there exists an $S$-scheme $S_0$ such that
\begin{displaymath}
\xymatrix{
\wtS\ar[dr]\ar[rr] & & S_0\ar[dl]\\
& S
}
\end{displaymath}
is a commutative diagram,  all morphisms are \'etale and $\left(\cC\times_SS_0\to S_0,\cC\times_SS_0\to X,\sigma_{0,\ast}\right)$ is a projective family of stable maps into $X$ (cfr. Proposition \ref{propA.2});
\item let $\wtgamma\colon\wtcC\to X\times_{\Z}\wtS$ be the morphism associated to the family $\left(\wtcC\to\wtS,\wtcC\to X,\wtsigma_{\ast}\right)$ (Remark \ref{remA.1} and cfr. proof of Proposition \ref{propA.2}), let $\pr_1\colon X\times_{\Z}\wtS\to X$ be a canonical projection. $q\colon\fE_{\wtcC}=\left(\pr_1\circ\wtgamma\right)^{\ast}\fE\to\fF$ is a quotient Higgs sheaf, whose underlying coherent sheaf is flat over $\wtS$.
\end{enumerate}
\noindent A morphism $\left(\cC_1\to S_1,\wtS_1\to S_1,\cC_1\to X,\sigma_{1,\ast},\fF_1,q_1\right)\to\left(\cC_2\to S_2,\wtS_2\to S_2,\cC_2\to X,\sigma_{2,\ast},\fF_2,q_2\right)$ is the triple $\left(f,\wtg,\psi\right)$ where:
\begin{enumerate}[A)]
\item $f\colon\left(\cC_1\to S_1,\cC_1\to X,\sigma_{1,\ast}\right)\to\left(\cC_2\to S_2,\cC_2\to X,\sigma_{2,\ast}\right)$ is a morphism in $\osfM_{g,m}(X,\beta)$;
\item $g\colon\cC_1\to\cC_2$ is the morphism included in the previous item;
\item by the Universal Property of fibre product, there exists a unique morphism ${\wtg\colon\wtcC_1\to\wtcC_2}$ such that the following diagram
\begin{displaymath}
\xymatrix{
\wtcC_1\ar[d]\ar[r]^{\wtg} & \wtcC_2\ar[d]\\
\wtS_1\ar[r] & \wtS_2
}
\end{displaymath}
is Cartesian in $\sfSch$;
\item $\psi$ is an isomorphism which makes commutative the following diagram
\begin{displaymath}
\xymatrix{
\left(\pr_1\circ\wtgamma_1\right)^{\ast}\fE\ar@{=}[d]\ar[rr]^(.625){q_1} & & \fF_1\ar[d]^{\psi}\\
\left(\pr_1\circ\wtgamma_2\circ\wtg\right)^{\ast}\fE\ar[rr]_(.625){\wtg^{\ast}\left(q_2\right)} & & \wtg^{\ast}\fF_2
},
\end{displaymath}
here: $\wtg^{\ast}\left(q_2\right)$ is the morphism of Higgs sheaves obtained by base change, $\wtgamma_i\colon\wtcC_i\to X\times_{\Z}\wtS_i$ is the morphism associated to these families as indicated above.
\end{enumerate}
\end{definition}
\noindent Let $\displaystyle\bfP\colon\sfQuot_{\fE/\osfM_{g,m}(X,\beta)}\to\osfM_{g,m}(X,\beta)$ be the functor which sends the objects $\left(\cC\to S,\wtS\to S,\cC\to X,\sigma_{\ast},\fF,q\right)$ to $\left(\cC\to S,\cC\to X,\sigma_{\ast}\right)$, and the morphisms $\left(f,\wtg,\psi\right)$ to $f$.
\begin{theorem}\label{th1.2}
$\bfP\colon\sfQuot_{\fE/\osfM_{g,m}(X,\beta)}\to\osfM_{g,m}(X,\beta)$ is a morphism of algebraic stacks\footnote{
The definitions of these objects given in \cite{TSP}, \cite{A:J} and \cite{C-M:W} are equivalent. We refer the reader to \cite[\href{https://stacks.math.columbia.edu/tag/076V}{tag 076V}]{TSP}, \cite[Definitions 4.1.3 and 4.1.6 and Theorem 4.2.1.(2)]{A:J} and \cite[Lemma B.16]{C-M:W}.
}.
\end{theorem}
\begin{prf}
Let $\left(\cC_1\to S_1,\cC_1\to X,\sigma_{1,\ast}\right)\to\left(\cC_2\to S_2,\cC_2\to X,\sigma_{2,\ast}\right)$ be a morphism in $\osfM_{g,m}(X,\beta)$, and let $\left(\cC_2\to S_2,\wtS_2\to S_2,\cC_2\to X,\sigma_{2,\ast},\fF_2,q_2\right)$ be an object of $\sfQuot_{\fE/\osfM_{g,m}(X,\beta)}$. By definition
\begin{displaymath}
\xymatrix{
\cC_1\ar[d]\ar[r]_g\ar@/^1pc/[rr] & \cC_2\ar[d]\ar[r] & X\\
S_1\ar[r]_f\ar@/^1pc/[u]^{\sigma_{i,1}} & S_2\ar@/_1pc/[u]_{\sigma_{i,2}}
}
\end{displaymath}
is a Cartesian diagram in in the category of algebraic spaces\footnote{
The definitions of these objects given in \cite{TSP}, \cite{A:J} and \cite{C-M:W} are equivalent. We refer the reader to \cite[\href{https://stacks.math.columbia.edu/tag/0GE0}{tag 0GE0}]{TSP}, \cite[Historical comments at the end of Section 4.1]{A:J} and \cite[Lemma B.14]{C-M:W}.
} $\sfAlgSp$. By hypothesis, $\wtS_2\to S_2$ is a morphism such that $\wtcC_2\to X\times_{\Z}\wtS_2$ is a family of stable maps into $X$. Let ${\wtS_1=\wtS_2\times_{S_2}S_1}$ and let $\wtcC_1\to X\times_{\Z}\wtS_1$ be the base change of $\cC_1$ over $\wtS_1$: this is family of stable maps into $X$ (cfr. proof of Proposition \ref{propA.2}). Let $\wtg\colon\wtcC_1\to\wtcC_2$ be the base change of $g$, let $\fF_1=\wtg^{\ast}\fF_2$ and ${q_1=\wtg^{\ast}q_2}$; $\left(f,\wtg,\psi\right)\colon\left(\cC_1\to S_1,\wtS_1\to S_1,\cC_1\to X,\sigma_{1,\ast},\fF_1,q_1\right)\to\left(\cC_2\to S_2,\wtS_2\to S_2,\cC_2\to X,\sigma_{2,\ast},\fF_2,q_2\right)$ is a morphism of objects in $\sfQuot_{\fE/\osfM_{g,m}(X,\beta)}$.\smallskip

\noindent Let $\left(\cC_1\to S_1,\cC_1\to X,\sigma_{1,\ast}\right)\to\left(\cC_2\to S_2,\cC_2\to X,\sigma_{2,\ast}\right)\to\left(\cC_3\to S_3,\cC_3\to X,\sigma_{3,\ast}\right)$ be morphisms of objects in $\osfM_{g,m}(X,\beta)$, let $\xymatrix{
\left(\cC_1\to S_1,\wtS_1\to S_1,\cC_1\to X,\sigma_{1,\ast},\fF_1,q_1\right)\ar[d]\\
\left(\cC_3\to S_3,\wtS_3\to S_3,\cC_3\to X,\sigma_{3,\ast},\fF_3,q_3\right) & \left(\cC_2\to S_2,\wtS_2\to S_2,\cC_2\to X,\sigma_{2,\ast},\fF_2,q_2\right)\ar[l]
}$ be morphisms of objects in $\sfQuot_{\fE/\osfM_{g,m}(X,\beta)}$. One has the following diagram in $\sfSch$
\begin{displaymath}
\xymatrix{
\wtcC_1\ar[d]\ar@/^1pc/[rr] & \wtcC_2\ar[d]\ar[r] & \wtC_3\ar[d]\\
\wtS_1\ar[r] & \wtS_2\ar[r] & \wtS_3
},
\end{displaymath}
since the right square is Cartesian in $\sfSch$, there exists a unique morphism $\wtcC_1\to\wtcC_2$ which makes commutative the whole of diagram; moreover $\wtcC_1\cong\wtcC_2\times_{\wtS_2}\wtS_1$. Repeating an analogous reasoning in $\sfAlgSp$, one has $\cC_1\cong\cC_2\times_{S_2}S_1$. Finally, $\fF_1$ is isomorphic to the pull-back of $\fF_2$ on $\wtcC_2\times_{\wtS_2}\wtS_1$, and $q_1$ is the composition of the inverse of this isomorphism with the pull-back of $q_2$. By definition, $\bfP\colon\sfQuot_{\fE/\osfM_{g,m}(X,\beta)}\to\osfM_{g,m}(X,\beta)$ is a \emph{prestack} (\cite[Definition 3.4.1]{A:J}).\smallskip

\noindent Let $\left(\cC\to S,\cC\to X,\sigma_{\ast}\right)$ be a family of stable maps into $X$ over a scheme $S$, by Proposition \ref{propA.2} there exists an \'etale morphism $\wtS\to S$ such that $\wtcC=\cC\times_S\wtS\to\wtS$ is projective; this last one induces a unique morphism $\bfQ\colon\sfC\to\osfM_{g,m}(X,\beta)$ of (algebraic) stacks; where $\sfC$ is the subcategory of $\sfSch$ such that the objects are schemes $\wtS\to S$ satisfying condition \ref{def1.1}.\ref{def1.1.b} and the morphisms are \'etale.\smallskip

\noindent Consider the fibre product $\sfQuot_{\fE/\osfM_{g,m}(X,\beta)}\times_{\osfM_{g,m}(X,\beta)}\sfC$ (\cite[Construction 3.4.33 and Theorem 3.4.35]{A:J}), the \emph{objects in} $\sfQuot_{\fE/\osfM_{g,m}(X,\beta)}$ ($\sfC$, respectively) \emph{over} $\left(\cC\to S,\cC\to X,\sigma_{\ast}\right)$ are all and the only $\left(\cC\to S,\wtS\to S,\cC\to X,\sigma_{\ast},\fF,q\right)$ ($S$-schemes $\wtS$ satisfying condition \ref{def1.1}.\ref{def1.1.b}, respectively)  By all this, the \emph{objects of} ${\sfQuot_{\fE/\osfM_{g,m}(X,\beta)}\times_{\osfM_{g,m}(X,\beta)}\sfC}$ \emph{over} $\left(\cC\to S,\wtS\to S,\cC\to X,\sigma_{\ast},\fF,q\right)$ and $\wtS\to S$ are all and only the objects of the stack associated to the scheme $\Quot_{\fE_{\wtcC}/\wtcC/\wtS}$ (see Proposition \ref{prop1.1} and Remark \ref{rem1.1}). By definition of $\sfQuot_{\fE/\osfM_{g,m}(X,\beta)}$, the \emph{morphisms between the objects of} ${\sfQuot_{\fE/\osfM_{g,m}(X,\beta)}\times_{\osfM_{g,m}(X,\beta)}\sfC}$ \emph{over} $\left(\cC\to S,\wtS\to S,\cC\to X,\sigma_{\ast},\fF,q\right)$ and $\wtS\to S$ are all and the only morphisms between objects of the stack associated to the scheme $\Quot_{\fE_{\wtcC}/\wtcC/\wtS}$.\smallskip

\noindent It turns out that $\sfQuot_{\fE/\osfM_{g,m}(X,\beta)}\times_{\osfM_{g,m}(X,\beta)}\sfC$ over $\left(\cC\to S,\wtS\to S,\cC\to X,\sigma_{\ast},\fF,q\right)$ and $\wtS\to S$ is isomorphic to the stack associated to the scheme $\Quot_{\fE_{\wtcC}/\wtcC/\wtS}$; by \cite[Corollaries 4.13 and 6.32]{C-M:W}, Theorems \ref{th1.1} and \ref{thA.1}, one has the claim.
\end{prf}

\section{Higgs sheaves, slope semistability and base change}

\noindent Let $X$ be a smooth, projective $\K$-variety, let $H$ be a polarization of $X$ and let $\fE=(\cE,\varphi)$ be a torsion-free Higgs sheaf on $X$, if not otherwise indicated. We consider the characteristic class
\begin{displaymath}
\Delta(\cE)=\frac{1}{2r}c_2\left(\cE\otimes\cE^{\vee}\right)=c_2(\cE)-\frac{r-1}{2r}c_1(\cE)^2\in\rA^2(X)\otimes_{\Z}\Q=\rA^2(X)_{\Q},
\end{displaymath}
called \emph{discriminant} of $\cE$; here $r=\rank(\cE)$.
\begin{remark}
$c_k(\cE)$ is the $k$-th Chern class of $\cE$\footnote{
I advise that, by \cite[Theorem 1.9 and Section 11]{M:YI}, one can apply the usual theory of Chern classes for locally free sheaves to coherent sheaves on $X$.
}, this is an element of $\rA_k(X)$, the Abelian group of $k$-cycles on $X$ modulo \emph{rational equivalence}. As it customs, $\rA^{n-k}(X)=\rA_k(X)$ by definition.\smallskip

\noindent Finally, it is defined a morphism $\displaystyle\int_X\colon\rA^n(X)\to\Z$ of Abelian groups (\cite[Definition 1.4]{F:W}).
\end{remark}
\begin{theorem}[{\cite[Theorem 7]{L:A}}]
Let $\fE=(\cE,\varphi)$ be a $H$-slope semistable Higgs sheaf on $X$. Then
\begin{displaymath}
\int_X\Delta(\cE)\cap H^{n-2}\geq0.
\end{displaymath}
\end{theorem}
\noindent About the ``extremal'' case, the following theorem holds.
\begin{theorem}[{see \cite[Theorems 1.2 and 1.3]{B:HR} and \cite[Proposition 3.2]{L:LG}}]\label{th2.1}
Let $\fE=(E,\varphi)$ be a Higgs bundle over $X$.
\begin{enumerate}[a)]
\item\label{th2.1.a} If $\fE$ is slope semistable with respect to some polarization $H$ and $\displaystyle\int_X\Delta(\cE)\cap H^{n-2}=0$. Then $\fE$ is curve semistable.
\item If $\fE$ is curve semistable. Then $\fE$ is slope semistable with respect to some polarization $H$.
\end{enumerate} 
\end{theorem}
\begin{remark}\label{rem2.2}
In statement \ref{th2.1.a} of the previous Theorem, one can write $\Delta(E)\equiv_{num}0$ (\emph{i.e.} $\Delta(E)$ is \emph{numerically equivalent} to $0$) instead of $\displaystyle\int_X\Delta(\cE)\cap H^{n-2}=0$. Indeed, by replacing $E$ with $E\otimes E^{\vee}$ if needed, one may assume $c_1(E)=0$. By hypothesis, this allows one to apply \cite[Corollary 6, see also note after Theorem 11]{L:A}, so that $c_k(E)\equiv_{num}0$ for all $k>0$, and then $\Delta(E)\equiv_{num}0$.
\end{remark}
\noindent If $\varphi$ vanishes then any curve semistable vector bundle $E$ is semistable with respect to some polarization and $\Delta(E)\equiv_{num}0$ (cfr. Remark \ref{rem0.2}.\ref{rem0.2.a}). Otherwise one posits the following conjecture.
\begin{conjecture}[\textbf{Bruzzo and Gra\~na Otero Conjecture}]\label{conj}
Let $\fE$ be a curve semistable Higgs bundle over $X$. Then $\fE$ is slope semistable with respect to some polarization $H$ and $\displaystyle{\int_X\Delta(E)\cap H^{n-2}=0}$.
\end{conjecture}
\begin{remark}\label{rem2.1}
Repeating the reasoning of Lanza and Lo Giudice contained in \cite[Subsection 3.1]{L:LG}, using \cite[Corollary 6]{L:A} instead of \cite[Lemma 3.7]{S:CT:1}, one has that this Conjecture holds if it holds for any smooth, projective surface.
\end{remark}
\noindent Here we prove the following theorem, which is the first main result of this paper.
\begin{theorem}\label{th2.2}
If Bruzzo and Gra\~na Otero Conjecture holds for any smooth, complex, projective surfaces. Then it holds for any smooth, projective surfaces defined over $\K$.
\end{theorem}
\noindent In order to prove the previous theorem, we premise the following lemmata.
\begin{lemma}[{\cite[Lemma 2.2.1]{AC:PhD}}]\label{lem2.1}
Let $\fE=(\cE,\varphi)$ be a Higgs sheaf on a projective variety $X$. Then there exist an algebraically closed subfield $\K_0$ of $\K$, a variety $X_{\K_0}$ over $\K_0$ and a Higgs sheaf $\fE_0$ over $X_{\K_0}$, such that $\K_0$ is isomorphic to a subfield of $\C$, the following diagram is Cartesian
\begin{displaymath}
\xymatrix{
X\ar[r]^{f_0}\ar[d] & X_{\K_0}\ar[d]\\
\spec(\K)\ar[r] & \spec\left(\K_0\right)
}
\end{displaymath}
and $\fE=f_0^{\ast}\fE_0$.
\end{lemma}
\begin{prf}
By definition $\displaystyle X=\proj\left(\frac{\K\left[x_0,\dotsc,x_N\right]}{J}\right)$ for some $N\in\N_{\geq1}$. To give a Higgs sheaf on $X$ is equivalent to give a triple $\left\{U_i,\lambda_{ij},\varphi_i\right\}_{i,j\in I}$ where
\begin{itemize}
\item $I$ is a finite set of indexes, because $X$ is quasi-compact as topological space;
\item $U_i$'s are open affine subsets of $X$ which cover it;
\item $\cE_i=\cE_{\vert U_i}$, $\lambda_{ij}\colon\cE_{i\vert U_{ij}}\to\cE_{j\vert U_{ji}}$'s are the transition functions of $\cE$, where $U_{ij}=U_i\cap U_j$ are non-empty, affine, open subsets of $X$;
\item $\varphi_i=\varphi\left(U_i\right)\colon\cE\left(U_i\right)\to\left(\cE\otimes_{\cO_X}\Omega^1_X\right)\left(U_i\right)=\cE\left(U_i\right)\otimes_{\cO_X\left(U_i\right)}\Omega^1_X\left(U_i\right)$, this last equality holds because $\cE$ and $\Omega^1_X$ are coherent $\cO_X$-modules and the $U_i$'s are open, affine subsets of $X$.
\end{itemize}
By \cite[Proposition II.5.2.(a)]{H:RC:1}, the morphisms $\lambda_{ij}$'s and $\varphi_i$'s correspond to morphisms $\widetilde{\lambda_{ij}}$'s and $\widetilde{\varphi_i}$'s of opportune coherent modules.\smallskip

\noindent Since all these are modules of finite type over $\K$, one can consider the set $S$ of all coefficients of the polynomials which generate the ideals of these modules. By Noetherianity of these modules, $S$ is a finite set and let $\alpha_1,\dotsc,\alpha_p,\tau_1,\dots,\tau_q$ its elements, where each $\alpha_i$ is algebraic on $\Q$ and the $\tau_j$'s are algebraically independent on $\Q$ (Definition \ref{defB.1}.\ref{defB.1.b}); if one of these types of elements does not occur one has either $p=0$ or $q=0$. Let $\K_0$ be the algebraically closure of the field generated by $S$ over $\Q$; $\K_0$ is a subfield of $\C$. Indeed, by definition
\begin{displaymath}
\K_0=\overline{Q\left(\Q(\alpha_1,\dotsc,\alpha_p)[\tau_1,\dotsc,\tau_q]\right)}^{alg}
\end{displaymath}
where $Q(\_)$ is the quotient field of $\_$. Let $t_1,\dots,t_q\in\C$ be transcendental numbers algebraically independent on $\Q$, then
\begin{displaymath}
\K_0\cong\overline{Q\left(\Q(\alpha_1,\dotsc,\alpha_p)[t_1,\dotsc,t_q]\right)}^{alg}\subsetneqq\C.
\end{displaymath}
Let $J_0$ be the ideal of $\K_0\left[x_0,\dotsc,x_N\right]$ generated by the generators of $J$ view as elements of this ring, let $\displaystyle X_{\K_0}=\proj\left(\frac{\K_0\left[x_0,\dotsc,x_N\right]}{J_0}\right)$. By construction one has the Cartesian diagram of the claim.\smallskip

\noindent Using the same reasoning, one defines a triple $\left\{U_i^0,\lambda_{ij}^0,\varphi_i^0\right\}_{i,j\in I}$ which determines a Higgs sheaf $\fE_0=\left(\cE_0,\varphi_0\right)$ on $X_{\K_0}$ such that $f_0^{\ast}\fE_0=\fE$.
\end{prf}
\begin{lemma}[{\textbf{``Lefschetz principle''-type for slope semistable Higgs sheaves}}, {\cite[Lemma 2.2.4]{AC:PhD}}]\label{lem2.2}
Let $\fE=(\cE,\varphi)$ be a torsion-free Higgs sheaf on $(X,H)$ and let $\F$ be a field extension of $\K$. Consider the following Cartesian diagram
\begin{displaymath}
\xymatrix{
X_{\F}\ar[r]^{f}\ar[d] & X\ar[d]\\
\spec(\F)\ar[r] & \spec(\K)
}
\end{displaymath}
then $f^{\ast}\fE=\left(f^{\ast}\cE,f^{\ast}\varphi\right)$ is slope semistable as Higgs sheaf if and only if $\fE=(\cE,\varphi)$ is slope semistable.
\end{lemma}
\begin{prf}
By base change, $f^{\ast}H$ is a polarization of $X_{\F}$ (see Remark \ref{remB.1}, Proposition \ref{propB.1} and \cite[\href{https://stacks.math.columbia.edu/tag/0D2P}{tag 0D2P}]{TSP}). As usual, for each coherent, torsion-free subsheaf $\cF$ of $f^{\ast}\fE$ of positive rank, one sets $\displaystyle\mu\left(\cF\right)=\frac{1}{\rank(\cF)}\int_{X_{\F}}c_1(\cF)\cap\left(f^{\ast}H\right)^{n-1}$ where $\dim\left(X_{\F}\right)=\dim(X)=n$ (see Lemma \ref{lemB.1}).\smallskip

\noindent If $\fE$ is slope unstable then there exists a torsion-free Higgs subsheaf $\fF=\left(\cF,\varphi_{\vert\cF}\right)$ of $\fE$ such that ${\mu(\fF)>\mu(\fE)}$, and without loss of generality, one may assume $\cF$ is reflexive. Since $f$ is a flat morphism then $f^{\ast}\cF$ is also reflexive (\cite[Proposition 1.8]{H:RC:2}) hence ${\mu\left(f^{\ast}\fF\right)>\mu\left(f^{\ast}\fE\right)}$ \emph{i.e.} $f^{\ast}\fE$ is unstable.\smallskip

\noindent If $f^{\ast}\fE$ is slope unstable, then there exists a torsion-free quotient Higgs sheaf $\fQ=\left(\cQ,\wtvarphi\right)$ such that $\mu(\fQ)<\mu\left(f^{\ast}\fE\right)$. From another point of view, $\fQ$ is a family of quotient Higgs sheaves of $\fE$ parameterised by $\spec(\F)$, it corresponds to a $\spec(\F)$-point of $\fQuot_{\fE/X/\K}$ (Theorem \ref{th1.1}) which belongs to an irreducible, closed subscheme $Z$. Since $\fQuot_{\fE/X/\K}$ is projective, $Z$ has a $\spec(\C)$-point which corresponds to a quotient Higgs sheaf $\fQ_0=\left(\cQ_0,\wtvarphi_0\right)$ of $\fE$ flat over $\spec(\K)$. By flatness, $\cQ_0$ is torsion-free and it has the same slope of $\cQ$, indeed, both corresponds to points of $Z$, hence they have the same Hilbert polynomial. Finally
\begin{displaymath}
\mu\left(\fQ_0\right)=\mu(\fQ)<\mu\left(f^{\ast}\fE\right)=\mu(\fE)
\end{displaymath}
which means that $\fE$ is a slope unstable Higgs sheaf.
\end{prf}
\begin{corollary}[{\textbf{``Lefschetz principle''-type for curve semistable Higgs sheaves}}]\label{cor2.1}
Without changing the hypothesis of the previous lemma, $f^{\ast}\fE=\left(f^{\ast}\cE,f^{\ast}\varphi\right)$ is a curve semistable Higgs sheaf if and only if $\fE=(\cE,\varphi)$ is curve semistable too.
\end{corollary}
\begin{prf}
Let $g\colon C\to X$ be a morphism, then $C_{\F}\equiv C\times_{\K}\spec(\F)$ is a curve (Lemma \ref{lemB.1}). By the previous Lemma, if $f^{\ast}\fE_{\vert C_{\F}}$ is slope semistable then $\fE_{\vert C}$ is slope semistable. Thus if $f^{\ast}\fE$ is curve semistable then $\fE$ is curve semistable as well. \emph{Vice versa}, let $\fE$ be curve semistable. Let ${g_0\colon C_0\to X_{\F}}$ be a morphism. If $\left(f\circ g_0\right)^{\ast}\fE=\fE_0$ is unstable, let $g\left(C_0\right)$ be the arithmetic genus of $C_0$. If $g\left(C_0\right)=0$ one chooses $3$ pairwise distinct closed points of $C_0$; if $g\left(C_o\right)=1$ one chooses $1$ closed point of $C_0$; otherwise no choice is necessary. In this way, $C_0\to X\times_{\C}\F$ is family of stable maps into $X$ parameterised by $\spec(\F)$. Let $q_0\colon\fE_0\to\fQ_0=\left(\cQ_0,\wtvarphi_0\right)$ be a quotient Higgs sheaf which destabilises $\fE_0$; by definition $\left(C_0\to\spec(\F),\spec(\F)\xrightarrow{\Id}\spec(\F),C_0\xrightarrow{f\circ g_0}X,\sigma_{\ast},\fQ_0,q_0\right)$ is an object of $\sfQuot_{\fE/\osfM_{g,m}(X,\beta)}$ (for some $\beta\in\rA^{n-1}(X)$) over $\spec(\F)$. Via the canonical morphism $\spec(\F)\to\spec(\C)$, there exists an object $\left(C\to\spec(\C),\spec(\C)\xrightarrow{\Id}\spec(\C),C\xrightarrow{g}X,\sigma_{\ast},\fQ,q\right)$ in $\sfQuot_{\fE/\osfM_{g,m}(X,\beta)}$. By the proof of Theorem \ref{th1.2}, $(\fQ,q)=\left(\cQ,\wtvarphi,q\right)$ determines an object over $\spec(\C)$ of the stack associated to the scheme $\fQuot_{\fE_0/C_0/\F}^{\deg\left(\cQ_0\right)\lambda+\rank\left(\cQ_0\right),\cL}$; where $\cL$ is an ample invertible sheaf on $C_0$. \emph{A priori}, $C$ is a connected, projective curve at worst nodal; let $\nu\colon\wtC\to C$ be its normalization, let $\wtC_1,\dotsc,\wtcC_{\ell}$ be its irreducible components and let $\nu_1=\nu_{\vert\wtcC_1},\dotsc,\nu_{\ell}=\nu_{\vert\wtcC_{\ell}}$. One has $\displaystyle\rA^1\left(\wtC\right)=\bigoplus_{k=1}^{\ell}\rA^1\left(\wtC_k\right)$, so
\begin{displaymath}
\deg\left(g^{\ast}\cE\right)=\sum_{k=1}^{\ell}\deg\left(\left(g\circ\nu_k\right)^{\ast}\cE\right)\,\text{\,and\,}\,\deg(\cQ)=\sum_{k=1}^{\ell}\deg\left(\nu_k^{\ast}\cQ\right).
\end{displaymath}
Thus
\begin{displaymath}
\frac{1}{\rank(\cQ)}\sum_{k=1}^{\ell}\deg\left(\nu_k^{\ast}\cQ\right)=\mu(\fQ)=\mu\left(\fQ_0\right)<\mu\left(\left(f\circ g_0\right)^{\ast}\fE\right)=\mu\left(g^{\ast}\fE\right)=\frac{1}{\rank(\cE)}\sum_{k=1}^{\ell}\deg\left(\left(g\circ\nu_k\right)^{\ast}\cE\right)
\end{displaymath}
and this is possible if and only if there exists $k_0\in\{1,\dotsc,\ell\}$ such that $\mu\left(\nu_{k_0}^{\ast}\fQ\right)<\mu\left(\left(g\circ\nu_{k_0}\right)\fE\right)$; in other words, $\fE$ is not curve semistable in contradiction with the assumptions. In order to avoid this absurd, $f^{\ast}\fE$ is curve semistable as well $\fE$.
\end{prf}
\noindent As application of Corollary \ref{cor2.1}, finally we are in position to prove Theorem \ref{th2.2}.\medskip

\noindent \textbf{Proof of Theorem \ref{th2.2}.} Let $X$ be a smooth, projective surface over $\K$, let $\fE=(E,\varphi)$ be a curve semistable Higgs bundle. Since $\displaystyle\Delta(E)=\frac{1}{2r}c_2\left(E\otimes E^{\vee}\right)$, one may replace $\fE$ by $\fE\otimes\fE^{\vee}$ hence to assume $c_1(E)=0$. So one has to prove only that $\displaystyle\int_Xc_2(E)\cap\_=0$. Repeating the proof of Lemma \ref{lem2.1}, one has a Cartesian diagram
\begin{displaymath}
\xymatrix{
X\ar[r]^{f_0}\ar[d] & X_{\K_0}\ar[d]\\
\spec(\K)\ar[r] & \spec\left(\K_0\right)
}
\end{displaymath}
where the notations are explained in that lemma; in particular, $\fE=f_0^{\ast}\fE_0$. By Corollary \ref{cor2.1}, $\fE_0$ is a curve semistable Higgs bundle. Without loss of generality, one may assume $\K=\K_0$ which is a subfield of $\C$, up to isomorphism. One changes the base of $X$ and has the following Cartesian diagram
\begin{displaymath}
\xymatrix{
X_{\C}=X\times_{\spec\left(\K\right)}\spec(\C)\ar[r]^(.75){f}\ar[d] & X\ar[d]\\
\spec(\C)\ar[r] & \spec\left(\K\right)
}.
\end{displaymath}
\noindent By Lemma \ref{lemB.1}, $X_{\C}$ is a smooth, complex, projective surface. Since $f$ is a flat morphism of schemes, one has the following diagram where the upper square is commutative (\cite[\href{https://stacks.math.columbia.edu/tag/0FVN}{tag 0FVN}]{TSP})
\begin{displaymath}
\xymatrix{
\Z\cong\rA^0(X)\ar[d]_{c_2(E)\cap\_}\ar[r]^(.45){\Id} & \rA^0\left(X_{\C}\right)\cong\Z\ar[d]^{c_2\left(f^{\ast}E\right)\cap\_}\\
\rA^2(X)\ar[r]_{f_0^{\ast}}\ar[d]_{\int_X\_} & \rA^2\left(X_{\C}\right)\ar[d]^{\int_{X_{\C}}\_} & .\\
\Z\ar@{=}[r] & \Z}
\end{displaymath}
Let $\displaystyle[Z]=\sum_{k=1}^mn_k\left[P_k\right]\in\rA^2(X)$, by definition
\begin{displaymath}
\int_X[Z]=\sum_{k=1}^mn_k\dim_{\K}\left(\kappa\left(P_k\right)\right)=\sum_{k=1}^mn_k\dim_{\K}\left(\rH^0\left(\left\{P_k\right\},\cO_{P_k}\right)\right);
\end{displaymath}
since $f$ is flat, by \cite[\href{https://stacks.math.columbia.edu/tag/02KH}{tag 02KH}]{TSP} one has
\begin{gather*}
\int_X[Z]=\sum_{k=1}^mn_k\dim_{\K}\left(\rH^0\left(\left\{P_k\right\},\cO_{P_k}\right)\right)=\sum_{k=1}^mn_k\dim_{\C}\left(\rH^0\left(\left\{P_k\right\},\cO_{P_k}\right)\right)\otimes_{\K}\C=\\
=\sum_{k=1}^mn_k\dim_{\C}\left(\rH^0\left(f^{-1}\left(P_k\right),\cO_{f^{-1}\left(P_k\right)}\right)\right)=\sum_{k=1}^mn_k\int_{X_{\C}}\left[f^{-1}\left(P_k\right)\right]=\int_{X_{\C}}f^{\ast}_0([Z]).
\end{gather*}
Thus, since $\displaystyle\int_{X_{\C}}c_2\left(f^{\ast}E\right)\cap\_=0$ then it has to be $\displaystyle\int_Xc_2(E)\cap\_=0$ and this concludes the proof.~\hfill{(Q.e.d.)}

\section{Numerically flat Higgs bundles}

\noindent Let $\fE=(E,\varphi)$ be a rank $r\geq2$ Higgs bundle over a smooth, projective $\K$-variety $X$, and let $s\in\{1,\dotsc,r-1\}$. Let $p_s\colon\Gr_s(E)\to X$ be the \emph{Grassmann bundle} parameterising rank $s$ locally free quotients of $E$ (see \cite{L:RK}). Consider the short exact sequence of vector bundles over $\Gr_s(E)$
\begin{displaymath}
\xymatrix{
0\ar[r] & S_{r-s,E}\ar[r]^(.55){\eta} & p_s^{\ast}E\ar[r]^(.45){\epsilon} & Q_{s,E}\ar[r] & 0
},
\end{displaymath}
where $S_{r-s,E}$ is the \emph{universal rank $r-s$ subbundle} and $Q_{s,E}$ is the \emph{universal rank $s$ quotient bundle of} $p_s^{\ast}E$, respectively. One defines the closed subschemes $\hGr_s(\fE)\subseteq\Gr_s(E)$ (the $s$\emph{-th Higgs-Grassmann schemes of} $\fE$) as the zero loci of the composite morphisms
\begin{displaymath}
\xymatrix{
(\epsilon\otimes\Id)\circ p_s^{\ast}\varphi\circ\eta\colon S_{r-s,E}\ar[r] & Q_{s,E}\otimes p_s^{\ast}\Omega_X^1
}.
\end{displaymath}
The restriction of the previous sequence to $\hGr_s(\fE)$ yields a universal short exact sequence
\begin{displaymath}
\xymatrix{
0\ar[r] & \mathfrak{S}_{r-s,\fE}\ar[r] & \rho_s^{\ast}\fE\ar[r] & \fQ_{s,\fE}\ar[r] & 0,
}
\end{displaymath}
where $\fQ_{s,\fE}=Q_{s,E|\hGr_s(\fE)}$ is equipped with the quotient Higgs field induced by $\rho_s^{\ast}\varphi$, and ${\rho_s=p_{s|\hGr_s(\fE)}}$. The scheme $\hGr_s(\fE)$ enjoys the usual \emph{Universal Property}: for a morphism of varieties $f\colon Y\to X$, the morphism $g\colon Y\to\Gr_s(E)$ given by a rank $s$ quotient $Q$ of $f^{\ast} E$ factors through $\hGr_s(\fE)$ if and only if $\varphi$ induces a Higgs field on $Q$.
\begin{definition}[{see \cite[Definition 2.3]{B:GO:2}}]\label{def3.1}
A Higgs bundle $\fE=(E,\varphi)$ of rank one is said to be \emph{Higgs numerically effective} (\emph{H-nef}, for short) if $E$ is \emph{numerically effective} in the usual sense\footnote{
Following \cite{L:RK}, a line bundle $L$ over a proper scheme $Z$ is \emph{numerically effective} (\emph{nef}, for short) if $\displaystyle\int_Cc_1(L)\geq0$ for every irreducible curve $C\subseteq Z$.
}. If $\rank(\fE)\geq 2$, one defines H-nefness inductively by requiring that
\begin{enumerate}[a)]
\item all Higgs bundles $\fQ_{s,\fE}$ are H-nef for all $s$, and
\item the determinant line bundle $\det(E)$ is numerically effective.
\end{enumerate}
If both $\fE$ and $\fE^{\vee}$ are Higgs numerically effective, $\fE$ is said to be \emph{Higgs numerically flat} (\emph{H-nflat}, for short).
\end{definition}
\begin{remarks}\label{rem3.1}
\,\begin{enumerate}[a)]
\item In the previous definition, the condition on the determinant cannot be omitted as \cite[Example 2.5]{B:GO:2} proves.
\item\label{rem3.1.b} Let $\fE=(E,\varphi)$ be a Higgs bundle over $X$. Where $\varphi=0$, $\fE$ is H-nef if and only if $E$ is nef as ordinary vector bundle. \hfill{$\Diamond$}
\end{enumerate}
\end{remarks}
\begin{proposition}\label{prop3.1}
A Higgs bundle $\fE=(E,\varphi)$ is H-nflat if and only if it is curve semistable and $\displaystyle\int_Cf^{\ast}c_1(E)=0$ for any $f\colon C\to X$ (\cite[Lemma 3.3.8]{AC:PhD}).
\end{proposition}
\noindent H-nflat Higgs bundles are related to Conjecture \ref{conj} by the following theorem.
\begin{theorem}[{cfr. \cite[Corollary 3.2]{B:B:G}}]\label{th3.1}
Bruzzo and Gra\~na Otero Conjecture holds for $X$ if and only if the Chern classes of any H-nflat Higgs bundles over $X$ are numerically trivial.
\end{theorem}
\begin{prf}
It is enough repeat the proof of \cite[Corollary 3.2]{B:B:G} and, taking into account Remark \ref{rem2.2}, one applies \cite[Corollary 6]{L:A} instead of \cite[Theorem 2]{S:CT:1}.
\end{prf}

\subsection{Semistability conditions for H-nflat Higgs bundles over smooth, complex, projective surfaces}

\noindent In order to prove Conjecture \ref{conj} on complex elliptic surfaces, henceforward we assume without loss of generality $\dim(X)=2$ and $\K=\C$, unless otherwise indicated.\medskip

\noindent For the aim of this paper, we prove the following lemma which works also for smooth projective varieties of any dimension, defined over an algebraically closed field of characteristic $0$.
\begin{lemma}[{cfr. \cite[Remark at page 17]{S:CT:3}}]\label{lem3.2}
Provided that H-nflat Higgs bundles over $X$ of rank at most $r-1$ have numerically trivial Chern classes. Then rank $r$ H-nflat Higgs bundles over $X$ are semistable, where $r\in\N_{\geq1}$.
\end{lemma}
\begin{prf}
By Remark \ref{rem0.1}, H-nflat Higgs line bundles are both stable and slope stable; hence let $r\in\N_{\geq2}$ for the rest of the proof.\smallskip

\noindent If $\fE$ is slope stable, there is nothing to prove by Remark \ref{rem0.1}. Otherwise, if $\fE$ is slope semistable, by \cite[Theorem 4.3.5]{AC:PhD} there exists a filtration
\begin{displaymath}
\uzero=\fE_0\subsetneqq\fE_1\subsetneqq\dotsc\subsetneqq\fE_{m-1}\subsetneqq\fE_m=\fE
\end{displaymath}
of Higgs subbundles such that each successive quotient $\fQ_k=\fE_k/\fE_{k-1}$ is locally free, (slope) stable and H-nflat. By the hypothesis, $\fQ_k$'s have the same normalised Hilbert polynomial; since $\fE$ is an iterated extension of these Higgs bundles, by \cite[Corollary 3.9]{HC:MG}, one has the claim.
\end{prf}

\subsection{On nilpotent, H-nflat Higgs bundles over smooth, complex, projective surfaces}

\noindent Let $\rM^{p,\cL}_X$ be the moduli space of Higgs sheaves on $X$, with fixed Hilbert polynomial $p\in\Q[\lambda]$ with respect to an ample invertible sheaf $\cL$ on $X$ and $\deg(p)=\dim(X)$ (see Appendix \ref{HiggsModSp}) and let
\begin{equation}\label{eq3.1}
\cS^{p,\cL}_X=\left\{[\fE]_{\sim}\in\rM^{p,\cL}_X\mid\fE\,\text{\,is JH-equivalent to an H-nflat Higgs bundle}\right\}
\end{equation}
\noindent the \emph{H-nflat locus}. This locus is well-defined whether H-nflat Higgs bundles are semistable; for example, the previous lemma gives a sufficient condition in order that this happens.
\begin{proposition}\label{prop3.2}
Provided that $\cS^{p,\cL}_X$ is well-defined and not empty, then it is stable under the action of $\C^{\times}$ on the moduli space $\rM^{p,\cL}_X$ of Higgs sheaves on $X$.
\end{proposition}
\begin{prf}
Let $f\colon C\to X$ be a morphism. This induces a morphism $f^{\sharp}\colon\cS^{p,\cL}_X\to\rM_C^{q,f^{\ast}\left(\cL\right)}$, where $g(C)$ is the genus of $C$ and $q\in\Q[\lambda]$ opportune. Trivially:
\begin{gather*}
\forall z\in\C^{\times},\,f^{\sharp}\left(z\cdot\left([\fE]_{\sim}\right)\right)=f^{\sharp}\left(z\cdot\left([(E,\varphi)]_{\sim}\right)\right)=f^{\sharp}\left([(E,z\varphi)]_{\sim}\right)=\left[\left(f^{\ast}E,f^{\ast}(z\varphi)\right)\right]_{\sim}=\\
=\left[\left(f^{\ast}E,zf^{\ast}\varphi\right)\right]_{\sim}=z\cdot\left(\left[\left(f^{\ast}E,f^{\ast}\varphi\right)\right]_{\sim}\right)=z\cdot\left(\left[f^{\ast}\fE\right]_{\sim}\right)=z\cdot\left(f^{\sharp}\left([\fE]_{\sim}\right)\right);
\end{gather*}
since there exists a H-nflat Higgs bundle $\fE_0$ over $X$ such that $[\fE]_{\sim}=\left[\fE_0\right]_{\sim}$ then
\begin{displaymath}
\forall z\in\C^{\times},\,f^{\sharp}\left(z\cdot\left(\left[\fE_0\right]_{\sim}\right)\right)=f^{\sharp}\left(z\cdot\left([\fE]_{\sim}\right)\right)=z\cdot\left(f^{\sharp}\left([\fE]_{\sim}\right)\right)=z\cdot\left(f^{\sharp}\left(\left[\fE_0\right]_{\sim}\right)\right),
\end{displaymath}
by Proposition \ref{prop3.1}, $f^{\sharp}\left(z\cdot\left(\left[\fE_0\right]_{\sim}\right)\right)$ is the JH-equivalence class of a H-nflat Higgs bundle over $C$ for any choice of $z$, $C$ and $f$. Thus $z\cdot\left([\fE]_{\sim}\right)$ is JH-equivalent to a H-nflat Higgs bundle over $X$ for any $z\in\C^{\times}$, \emph{i.e.} the claim holds.
\end{prf}
\begin{corollary}\label{cor3.1}
Without change the previous notations, the limit of a H-nflat Higgs bundle under the $\C^{\times}$-action on $\rM^{p,\cL}_X$ is a nilpotent, H-nflat Higgs bundle.
\end{corollary}
\begin{prf}
By Remark \ref{remC.1}, limits of Higgs bundles in $\cS^{p,\cL}_X$ are nilpotent Higgs bundles. Repeating the proof of the previous proposition, one has
\begin{displaymath}
f^{\sharp}\left(\lim_{z\to0}z\cdot\left(\left[\fE\right]_{\sim}\right)\right)=\lim_{z\to0}f^{\sharp}\left(z\cdot\left([\fE]_{\sim}\right)\right)=\lim_{z\to0}z\cdot\left(f^{\sharp}\left(\left[\fE\right]_{\sim}\right)\right),
\end{displaymath}
where $\fE$ is a H-nflat Higgs bundle over $X$, $f\colon C\to X$ is a morphism, $f^{\sharp}\colon\cS^{p,\cL}_X\to\rM_C^{q,f^{\ast}\left(\cL\right)}$, where $g(C)$ is the genus of $C$ and $q\in\Q[\lambda]$ opportune. $z\cdot\left(f^{\sharp}\left(\left[\fE\right]_{\sim}\right)\right)$ is the JH-equivalence class of a semistable Higgs bundle over $C$ with $c_1(E)=0\in\rA^1(C)$; since this holds for any $z$, $f$ and $C$, then $\displaystyle\lim_{z\to0}z\cdot\left(\left[\fE\right]_{\sim}\right)$ is the JH-equivalence class of H-nflat Higgs bundle, \emph{i.e.} the claim holds.
\end{prf}
\noindent Hence Conjecture \ref{conj} is simplified as it follows (see also Theorem \ref{th3.1}).
\begin{lemma}\label{lem3.1}
Provided that $\cS^{p,\cL}_X$ is well-defined and nilpotent, H-nflat Higgs bundles over $X$ have numerically trivial Chern classes, then Conjecture \ref{conj} holds for $X$.
\end{lemma}
\begin{prf}
Trivial, because $[\fE]_{\sim}\in\cS^{p,\cL}_X$ and its ``limit'' $\left[\fE_0\right]_{\sim}$ have the Chern classes.
\end{prf}

\section{H-nflat Higgs bundles over elliptic surfaces}

\subsection{Preliminaries about elliptic surfaces}

\noindent We fix the following definitions. Let $\pi\colon X\to B$ be a proper, flat, surjective morphism; where $B$ and $X$ are smooth, complex, projective varieties with $\dim(B)=1$ and $\dim(X)=2$.
\begin{definition}
If the generic scheme-theoretic fibre of $\pi$ is a smooth, connected curve of genus $1$ then $\pi$ is called \emph{elliptic fibration} (\emph{over} $B$) and $(X,\pi)$ is called \emph{elliptic surface}.
\end{definition}
\begin{definition}\label{def4.1}
An elliptic surface $(X,\pi)$ is \emph{Jacobian} if $\pi$ has a section $\sigma\colon B\to X$, \emph{i.e.} $\pi\circ\sigma=\Id_B$.
\end{definition}
\noindent Where there is no confusion, we shall write $X$ to indicate an elliptic surface. Since we are assuming $\K=\C$, we shall apply \emph{GAGA Principle} where it holds, tacitly.
\begin{lemma}[{cfr. \cite[Subsubsection 2.3.1]{B:P}}]\label{lem4.1}
The \emph{cotangent exact sequence}
\begin{displaymath}
\xymatrix{
\pi^{\ast}\Omega^1_B\ar[r]^{\pi^{\ast}} & \Omega^1_X\ar[r] & \Omega^1_{X/B}\ar[r] & 0
}
\end{displaymath}
is left exact.
\end{lemma}
\noindent The left exactness follows from the following known result, for which we are not able to indicate a reference.
\begin{proposition}\label{prop4.2}
$\Omega^1_{X/B}$ is a torsion-free, coherent sheaf on $X$ of generic rank $1$.
\end{proposition}
\begin{prf}
Since $\pi$ is a morphism of finite type between Noetherian schemes, reasoning via opportune affine open coverings of $X$ and $B$, one proves the coherence of $\Omega^1_{X/B}$. Let $X_{\eta}$ be the generic scheme-theoretic fibre of $\pi$. Easily $\left(\Omega^1_{X/B}\right)_{\vert X_{\eta}}\cong\Omega^1_{X_{\eta}}=\cO_{X_{\eta}}$, \emph{i.e.} $\Omega^1_{X/B}$ has generic rank $1$. By Noetherianity and irreducibility of $X$, $\Omega^1_{X/B}$ is torsion-free.
\end{prf}
\noindent \textbf{Proof of Lemma \ref{lem4.1}.} By generic rank counting, $\ker\left(\pi^{\ast}\right)$ has generic rank $1$, hence it is torsion-free. Since $\pi^{\ast}\Omega^1_B$ is reflexive, by \cite[Proposition 1.1]{H:RC:2}, the kernel of
\begin{displaymath}
\xymatrix{
\pi^{\ast}\Omega^1_B\ar[r] & \ker\left(\pi^{\ast}\right)\ar[r] & 0
}
\end{displaymath}
is a reflexive sheaf of generic rank $0$, \emph{i.e.} it is $\uzero_X$. By all this $\pi^{\ast}\Omega^1_B\cong\ker\left(\pi^{\ast}\right)$ and the cotangent exact sequence is also left exact, equivalently. \hfill{(Q.e.d.)}\medskip

\noindent Henceforth, let $X_{\eta}$ be the generic scheme-theoretic fibre of $\pi$. We recall that for each $r\in\N_{\geq1}$ there exists a unique up to isomorphisms, indecomposable, rank $r$ vector bundle $\I_r$ over any smooth, complex, projective curve $C$ of genus $1$ such that $\rh^0(C,\I_r)=1$ (\cite[Lemma 15 and Theorem 5.(i)]{A:MF}).
\begin{remark}\label{rem4.1}
Since $X_{\eta}$ is \emph{geometrically reduced} (see proof of Lemma \ref{lemB.1}), by \cite[\href{https://stacks.math.columbia.edu/tag/0578}{tag 0578}]{TSP}, there exists an open subset $U$ of $B$ such that $\pi$ has geometrically reduced fibres over $U$. In particular, there exist closed, reduced fibres of $\pi$.\medskip

\noindent If $X$ is Jacobian, then the fibres of $\pi$ are reduced (see \cite[Section II.3]{M:R} or \cite[Section 7.4]{F:R}).
\end{remark}
\begin{definition}
An elliptic surface $\pi\colon X\to B$ is \emph{isotrivial} if there exists a base change $\beta\colon C\to B$ such that $X\times_BC$ is birationally equivalent to $\Gamma\times_{\C}C$, where $\Gamma$ is an elliptic curve.
\end{definition}
\begin{remark}
Elliptic surfaces birationally equivalent are isomorphic (\cite[Proposition II.1.2]{M:R}).
\end{remark}
\noindent Elliptic surfaces $X$ with $\chi\left(\cO_X\right)>0$ are not isotrivial; in order to prove this fact, explicitly, we premise the following known statement for which we are not able to indicate some reference. These will be applied also in the next.
\begin{proposition}\label{prop4.1}
A morphism $f\colon Y\to Z$ of projective schemes over a scheme $S$ is projective\footnote{
For clarity, $f$ is \emph{projective} if there exist $M\in\N_{\geq1}$ and a closed immersion $i\colon Y\hookrightarrow\bP^M_Z$ such that $f$ factorises as $Y\stackrel{i}{\hookrightarrow}\bP^M_Z\to Z$ (\cite[Definition at page 103]{H:RC:1}).
}.
\end{proposition}
\begin{prf}
By hypothesis and \cite[Proposition 7.4.13]{B:S}, the \emph{graph morphism} $\Gamma_f\colon Y\to Y\times_SZ$ is a closed immersion. Since there exist $M\in\N_{\geq1}$ and a closed immersion $j\colon Y\hookrightarrow\bP^M_S$, one has a closed immersion $i\colon Y\stackrel{\Gamma_f}{\hookrightarrow}Y\times_SZ\stackrel{\widetilde{j}}{\hookrightarrow}\bP^M_Z$ via base change. Applying \emph{\emph{Universal Property} of fibre product of schemes}, $f$ factorises as $Y\stackrel{\Gamma_f}{\hookrightarrow}Y\times_SZ\stackrel{\widetilde{j}}{\hookrightarrow}\bP^M_Z\to Z$ and one has the claim.
\end{prf}
\begin{lemma}
An elliptic surface $X$ with $\chi\left(\cO_X\right)>0$ is not isotrivial.
\end{lemma}
\begin{prf}
Indeed, this conditions implies that $\bL_X=\left(R^1\pi_{\ast}\cO_X\right)^{\vee}$ is not the trivial line bundle (by \cite[Corollary III.12.9]{H:RC:1} and \cite[Corollaries 7.6 and 7.17]{F:R}). Given a base change $\beta\colon C\to B$, from flatness of $\pi$ follows that $\cO_X$ is a flat coherent sheaf on $B$; by the previous proposition and \cite[Theorem 25.1.5]{FOAG} $\bL_X\cong\bL_{X\times_BC}$. By all this, $X\times_BC$ is not isomorphic to $\Gamma\times_{\spec(C)}C$, where $\Gamma$ is an elliptic curve; because this last eventuality is equivalent to triviality of $\bL_{X\times_BC}$ (\cite[Lemma III.1.4]{M:R}).
\end{prf}
\noindent From now on, $X$ is a non-isotrivial elliptic surface, unless otherwise indicated. Under this assumption, the following lemma holds.
\begin{lemma}[{cfr. \cite[Subsubsection 2.3.2]{B:P}}]\label{lem4.2}
Let $X_b$ be a fibre of $X$, with $b\in B(\C)$ \emph{general}. Then $\Omega^1_{X\vert X_b}\cong\I_2$.
\end{lemma}
\begin{prf}
Let $U$ be a codimension $2$ open subset of $X$ on which $\Omega^1_{X/B}$ is locally free, one has the following extension of line bundles (Lemma \ref{lem4.1} and Proposition \ref{prop4.2})
\begin{displaymath}
\xymatrix{
0\ar[r] & \left(\pi^{\ast}\Omega^1_B\right)_{\vert U}\ar[r] & \Omega^1_{X\vert U}\ar[r] & \Omega^1_{X/B\vert U}\ar[r] & 0
}.
\end{displaymath}
For any $b\in B(\C)$ such that $X_b\subseteq U$ and $X_b$ is smooth of genus $1$ (\cite[Corollary III.9.10]{H:RC:1}), the previous sequence becomes
\begin{displaymath}
\xymatrix{
0\ar[r] & \cO_{X_b}\ar[r] & \Omega^1_{X\vert X_b}\ar[r] & \cO_{X_b}\ar[r] & 0
}.
\end{displaymath}
Since $X$ is non-isotrivial, by \cite[Lemma 2]{V:R} $\pi_{\ast}\Omega^1_X\cong\cO_B$, this is equivalent to $\rh^0\left(\Omega^1_{X_b}\right)=1$ by \cite[Corollary III.9.6]{H:RC:1}. Thus $\Omega^1_{X_b}\cong\I_2$.
\end{prf}

\subsection{On generally polystable Higgs bundles and vertical Higgs bundles over non-isotrivial elliptic surfaces}

\noindent Henceforth, $X_b$ is the scheme-theoretic fibre of $\pi$ over $b\in B(\C)$.
\begin{definition}
Let $\fE=(E,\varphi)$ be a Higgs bundle over $X$. It is \emph{generally} (\emph{semi})\emph{stable} (respectively, \emph{generally polystable}) if $\fE_{\vert X_b}$ is (semi)stable (respectively, polystable) for some $b\in B(\C)$ general. It is \emph{fibrewise} (\emph{semi})\emph{stable} (respectively, \emph{fibrewise polystable}) if $\fE_{\vert X_b}$ is (semi)stable (respectively, polystable) for any $b\in B(\C)$ general.
\end{definition}
\begin{lemma}[{cfr. \cite[Proposition 4.2]{B:P}}]
Let $\fE=(E,\varphi)$ be a rank $r\geq2$, generally semistable Higgs bundle over $X$. If $\displaystyle\gcd\left(\int_{X_b}c_1\left(E_{\vert X_b}\right),r\right)>1$ for some $b\in B(\C)$ general, then the natural morphism
\begin{displaymath}
\pi_{\ast}\left(\Id\otimes\pi^{\ast}\right)\colon\pi_{\ast}\End(E)\otimes_{\cO_X}\Omega^1_B\cong\pi_{\ast}\left(\End(E)\otimes_{\cO_X}\pi^{\ast}\Omega^1_B\right)\to\pi_{\ast}\left(\End(E)\otimes_{\cO_X}\Omega^1_X\right)
\end{displaymath}
is a monomorphism of vector bundles over $B$. Moreover, if $\fE$ is generally polystable then $\pi_{\ast}\left(\Id\otimes\pi^{\ast}\right)$ is an isomorphism of rank $r$ vector bundles over $B$.
\end{lemma}
\begin{prf}
we begin the proof by assuming that $\fE$ is generally polystable. Under this hypothesis, by \cite[Theorem 4.14]{F:GP:N} $\fE_{\vert X_b}$ is a polystable not stable Higgs bundle over $X_b$\footnote{
Since $b$ is a general closed point of $B$, without loss of generality, one may assume that $X_b$ is an irreducible, smooth, projective curve of genus $1$.
}, hence it is the direct sum of $r$ Higgs line bundles $\left(L_i,\psi_i\right)$ over $X_b$ of degree $\displaystyle\int_{X_b}c_1\left(E_{\vert X_b}\right)$.\smallskip

\noindent On another hand, since $\pi$ is flat then its fibres have all dimension $1$ (\cite[Corollary III.9.6]{H:RC:1}) hence $\pi_{\ast}\End(E)$ is a vector bundle over $B$ (\cite[Corollaries 1.4 and 1.7]{H:RC:2}).\smallskip

\noindent By flatness of $\pi$ follows also that $\End(\cE)$ is a flat coherent sheaf on $B$; by Proposition \ref{prop4.1} and by \cite[Corollary III.12.9]{H:RC:1}, it turns out that $\rh^0\left(X_b,\End\left(\cE_{\vert X_b}\right)\right)$ is constant and positive for any $b$. To be precise
\begin{displaymath}
H^0\left(X_b,\End\left(\cE_{\vert X_b}\right)\right)\cong H^0\left(X_b,\bigoplus_{i,j=1}^r\cL_i\otimes_{\cO_X}\cL_j^{\vee}\right)\cong H^0\left(X_b,\cO_{X_b}^{\oplus r}\right)
\end{displaymath}
hence for any $b\in B,\,\rh^0\left(X_b,\End\left(\cE_{\vert X_b}\right)\right)=r$ \emph{i.e.} $\rank(\pi_{\ast}\End(E))=r$. Moreover, repeating the same reasoning and by Lemma \ref{lem4.2}
\begin{displaymath}
H^0\left(X_b,\left(\End(\cE)\otimes\Omega^1_X\right)_{\vert X_b}\right)\cong H^0\left(X_b,\bigoplus_{i,j=1}^r\cL_i\otimes_{\cO_X}\cL_j^{\vee}\otimes\I_2\right)\cong H^0\left(X_b,\I_2^{\oplus r}\right)
\end{displaymath}
hence for any $b\in B,\,\rh^0\left(X_b,\left(\End(\cE)\otimes_{\cO_X}\Omega^1_X\right)_{\vert X_b}\right)=r$ \emph{i.e.} $\rank\left(\pi_{\ast}\left(\End(E)\otimes_{\cO_X}\Omega^1_X\right)\right)=r$.\smallskip

\noindent From all this, $\pi_{\ast}\left(\Id\otimes\pi^{\ast}\right)$ is a morphism of rank $r$ vector bundles over $B$. Finally, by exactness of the cotangent exact sequence (Lemma \ref{lem4.1}) one has the left exact sequences over $B$
\begin{displaymath}
\xymatrix{
0\ar[r] & \End(E)\otimes_{\cO_X}\pi^{\ast}\Omega^1_B\ar[r]^{\Id\otimes\pi^{\ast}} & \End(E)\otimes_{\cO_X}\Omega^1_X\ar[r] & \End(E)\otimes_{\cO_X}\Omega^1_{X/B}\ar[r] & 0,\\
0\ar[r] & \pi_{\ast}\left(\End(E)\otimes_{\cO_X}\pi^{\ast}\Omega^1_B\right)\ar[r]_(.525){\pi_{\ast}\left(\Id\otimes\pi^{\ast}\right)} & \pi_{\ast}\left(\End(E)\otimes_{\cO_X}\Omega^1_X\right)\ar[r] & \pi_{\ast}\left(\End(E)\otimes_{\cO_X}\Omega^1_{X/B}\right);
}
\end{displaymath}
by \emph{Projection Formula} $\pi_{\ast}\End(E)\otimes_{\cO_X}\Omega^1_B\cong\pi_{\ast}\left(\End(E)\otimes_{\cO_X}\pi^{\ast}\Omega^1_B\right)$, thus $\pi_{\ast}\left(\Id\otimes\pi^{\ast}\right)$ is a monomorphism of locally free sheaves of the same rank hence it is an isomorphism.\smallskip

\noindent More in general, let $\fE$ be generally semistable. Repeating the previous reasoning, one has that $\pi_{\ast}\End(E)$ is a vector bundle over $B$ and $\pi_{\ast}\left(\Id\otimes\pi^{\ast}\right)$ is a monomorphism of vector bundles over $B$.
\end{prf}
\noindent The previous lemma justifies the introduction of the following definition.
\begin{definition}[{see \cite{B:P}}]
Let $\fE=(E,\varphi)$ be a Higgs bundle over $X$. If $\varphi$ factorises through $E\otimes_{\cO_X}\pi^{\ast}\Omega^1_B$ then it is called \emph{vertical} (Higgs field), and $\fE$ is said \emph{vertical} (Higgs bundle).
\end{definition}
\begin{corollary}[{cfr. \cite[Corollary 4.3]{B:P}}]
Let $\fE=(E,\varphi)$ be a rank $r\geq2$, generally polystable Higgs bundle over $X$. If $\displaystyle\gcd\left(\int_{X_b}c_1\left(E_{\vert X_b}\right),r\right)>1$ for some $b\in B(\C)$ general, then $\varphi$ and $\varphi^{\vee}$ are both vertical.
\end{corollary}
\begin{prf}
To give a Higgs field on $E$ is equivalent to choose a global section of $\End(E)\otimes_{\cO_X}\Omega^1_X$; by construction this is a global section of $\pi_{\ast}\left(\End(E)\otimes_{\cO_X}\Omega^1_X\right)$. By the previous lemma, this corresponds to a global section of $\pi_{\ast}\left(\End(E)\otimes_{\cO_X}\pi^{\ast}\Omega^1_B\right)$ which is a global section of ${\End(E)\otimes_{\cO_X}\pi^{\ast}\Omega^1_B}$.\smallskip

\noindent To be clear, the dual of $\varphi$ is the opposite of the following composition
\begin{displaymath}
E^{\vee}\xrightarrow{\cong}E^{\vee}\otimes_{\cO_X}\cO_X\xrightarrow{\Id\otimes\Tr^{\vee}}E^{\vee}\otimes_{\cO_X}\left(\pi^{\ast}\Omega^1_B\right)^{\vee}\otimes_{\cO_X}\pi^{\ast}\Omega^1_B\xrightarrow{\varphi^{\vee}\otimes\Id}E^{\vee}\otimes_{\cO_X}\pi^{\ast}\Omega^1_B
\end{displaymath}
which is a vertical Higgs field, under the present hypothesis.
\end{prf}
\begin{lemma}
Let $\fE=(E,\varphi)$ be a rank $r\geq2$, generally polystable Higgs bundle over $X$. If $\displaystyle\gcd\left(\int_{X_b}c_1\left(E_{\vert X_b}\right),r\right)>1$ for some $b\in B(\C)$ general, then $\End(\fE)$ is not a slope stable Higgs bundle with respect to all polarizations of $X$.
\end{lemma}
\begin{prf}
For sake of simplicity, we set $\fE_0=\left(E_0,\varphi_0\right)=(\End(E),\End(\varphi))$. By hypothesis, $\fE_{\vert X_b}$ is a polystable Higgs bundle over $X_b$, by \cite[Lemma 4.4]{IB:GS} $\fE_{0,\vert X_b}$ is polystable too. Since $c_1\left(E_0\right)=0$, by the previous corollary $\varphi_0$ is a vertical Higgs bundle.\smallskip

\noindent Applying the \emph{Adjunction Formula} to the functors $\pi^{\ast}$ and $\pi_{\ast}$, the identity morphism of $\pi_{\ast}E_0$ corresponds to a canonical morphism $\gamma\colon\pi^{\ast}\left(\pi_{\ast}E_0\right)\to E_0$. Consider the following commutative diagram
\begin{displaymath}
\xymatrix{
\pi^{\ast}\pi_{\ast}E_0\ar[dr]_{\varphi_{\gamma}}\ar[r]^{\gamma} & E_0\ar[d]^{\varphi_0}\\
 & E_0\otimes\pi^{\ast}\Omega^1_B
};
\end{displaymath}
again, by the \emph{Adjunction Formula}, $\varphi_{\gamma}$ corresponds to a morphism $\widehat{\varphi}_{\gamma}\colon\pi_{\ast}E_0\to\pi_{\ast}\left(E_0\otimes_{\cO_X}\pi^{\ast}\Omega^1_B\right)$ and, by the \emph{Projection Formula}, this corresponds to a morphism $\wtvarphi_{\gamma}\colon\pi_{\ast}E_0\to\pi_{\ast}E_0\otimes_{\cO_X}\Omega^1_B$. By all this, one has the following diagram
\begin{displaymath}
\xymatrix{
\pi^{\ast}\pi_{\ast}E_0\ar[d]_{\pi^{\ast}\wtvarphi_{\gamma}}\ar[dr]^{\varphi_{\gamma}}\ar[r]^{\gamma} & E_0\ar[d]^{\varphi_0}\\
\pi^{\ast}\pi_{\ast}E_0\otimes_{\cO_X}\pi^{\ast}\Omega^1_B\ar[r]_(.55){\gamma\otimes\Id} & E_0\otimes_{\cO_X}\pi^{\ast}\Omega^1_B
}
\end{displaymath}
Explicitly, for any $V\subseteq B$, $U\subseteq X$ with $U\subseteq\pi^{-1}(V)$ opens, $\gamma$ is induced by
\begin{displaymath}
a\otimes b\in\cE_0\left(f^{-1}(V)\right)\otimes_{\cO_B(V)}\cO_X\left(f^{-1}(V)\right)\mapsto(ba)_{\vert U}\in\cE_0(U)
\end{displaymath}
(see \cite[at page 271]{B:S}) and so $\varphi_{\gamma}$ is induced by
\begin{displaymath}
a\otimes b\in\cE_0\left(f^{-1}(V)\right)\otimes_{\cO_B(V)}\cO_X\left(f^{-1}(V)\right)\mapsto\left(b\varphi_0\left(f^{-1}(V)\right)(a)\right)_{\vert U}\in\left(\cE_0\otimes_{\cO_X}\pi^{\ast}\omega^1_B\right)(U)
\end{displaymath}
hence $\wtvarphi_{\gamma}$ is defined as
\begin{displaymath}
\wtvarphi_{\gamma}(V)\colon a\in\cE_0\left(f^{-1}(V)\right)\to\varphi_{\gamma}\left(f^{-1}(V)\right)(a)\in\left(\cE_0\otimes_{\cO_X}\pi^{\ast}\omega^1_B\right)\left(f^{-1}(V)\right).
\end{displaymath}
Finally, $\pi^{\ast}\wtvarphi_{\gamma}$ is induced by
\begin{displaymath}
a\otimes b\in\cE_0\left(f^{-1}(V)\right)\otimes_{\cO_B(V)}\cO_X\left(f^{-1}(V)\right)\to\left(b\varphi_0\left(f^{-1}(V)\right)(a)\otimes1\right)_{\vert U}\in\left(\cE_0\otimes_{\cO_X}\pi^{\ast}\omega^1_B\right)(U)\otimes_{\cO_X}\cO_X(U)
\end{displaymath}
and by construction $\varphi_{\gamma}=(\gamma\otimes\Id)\circ\pi^{\ast}\wtvarphi_{can}$, \emph{i.e.} $\gamma\colon\pi^{\ast}\pi_{\ast}E_0\to E_0$ is a morphism of Higgs bundles.\smallskip

\noindent Since $\gamma$ is not the zero morphism, its image is a degree $0$, torsion-free, Higgs subsheaf of $\fE_0$ hence it cannot be slope stable.
\end{prf}
\begin{remark}
The previous lemma works also if $\End(\fE)$ is a vertical Higgs bundle.
\end{remark}
\noindent Last lemma allows us to generalise partially \cite[Theorem 4.14]{F:GP:N} to non-isotrivial elliptic surfaces.
\begin{theorem}
It does not exist Higgs bundles $\fE=(E,\varphi)$ over $X$ such that $\fE$ is slope stable with respect to some polarization $H$, $\fE_{\vert X_b}$ is polystable and $\displaystyle\gcd\left(\int_{X_b}c_1\left(E_{\vert X_b}\right),r\right)>1$, for some $b\in B(\C)$ general.
\end{theorem}
\begin{prf}
Indeed, if there exists a such Higgs bundle $\fE$ then $\End(\fE)$ is slope stable with respect to $H$ (\cite[Proposition 4.5]{IB:GS}); but this is impossible by the previous lemma.
\end{prf}
\begin{corollary}
Let $\fE$ be a slope stable, H-nflat Higgs bundle over $X$. Then $\fE_{\vert X_b}$ is semistable but not polystable for any $b\in B(\C)$ general.
\end{corollary}
\noindent The non-isotrivial hypothesis is necessary as the following example proves.
\begin{example}\label{ex4.1}
Let $B$ be a smooth, complex, projective curve of genus $1$, let ${P\in B(\C)}$ general and let $E$ be the non-split extension
\begin{displaymath}
\xymatrix{
0\ar[r] & \cO_B\ar[r] & E\ar[r] & \cO_B(P)\ar[r] & 0
}
\end{displaymath}
(see \cite[Theorem 2.6.(ii).(c)]{F:R}). Since there exists a degree $1$ line bundle $L$ over $B$ such that $E\cong\I_2\otimes_{\cO_B}L$ (\cite[Theorem 10]{A:MF}) then $E$ is stable by \cite[Lemma 30]{T:LW}. Moreover, this is equivalent to the stability of $\End(E)=E\otimes E^{\vee}$ (\cite[Proposition 4.5]{IB:GS}), and since $c_1(\End(E))=0$ then $\End(E)$ is nflat (Proposition \ref{prop3.1}).\smallskip

\noindent Let $C_g$ be a smooth, complex, projective curve of genus $g\in\N_{\geq0}$, let $X=B\times_{\C}C_g$ and let $\pr\colon X\to B$ be the canonical projection onto $B$. Applying \cite[Lemma 3.4]{B:B:G}, one has that ${F=\pr^{\ast}\End(E)\cong\pr^{\ast}E\otimes_{\cO_X}\pr^{\ast}E}$ is a nflat vector bundle of rank $4$.\smallskip

\noindent On the other hand, let $H$ be a polarization of $X$, by \emph{Bertini's Theorem} (\cite[Theorem II.8.18 and Remark III.7.9.1]{H:RC:1}), there exist $m\gg1$ and $Y\in\lvert mH\rvert$ such that $Y$ is a smooth, projective curve; by \cite[Theorem 6.1]{MVB:RA}, $\pr^{\ast}E_{\vert Y}$ is semistable. Let
\begin{displaymath}
\xymatrix{
q\colon Y\ar@{^{(}->}[r] & X\ar[r]^{\pr} & B,
}
\end{displaymath}
it is a finite morphism and $\pr^{\ast}E_{\vert Y}=q^{\ast}E$ trivially, \emph{i.e.} $q^{\ast}E$ fits in the following short exact sequence
\begin{displaymath}
\xymatrix{
0\ar[r] & q^{\ast}\cO_B=\cO_X\ar[r] & q^{\ast}E\ar[r] & q^{\ast}\cO_B(P)\ar[r] & 0.
}
\end{displaymath}
By hypothesis, $\displaystyle\int_Yc_1\left(q^{\ast}E\right)=\deg(f)\int_Bc_1(E)=\deg(f)\geq1$ and, by \cite[Theorem 1]{K:RG}, $q^{\ast}E$ is polystable. If $q^{\ast}E$ is not stable then there exist $L_1$ and $L_2$ line bundles over $Y$ such that ${q^{\ast}E=L_1\oplus L_2}$ and $\displaystyle\int_Yc_1\left(L_1\right)=\int_Yc_1\left(L_2\right)=\deg(f)$; but it has to be $L_k=\cO_Y$ for some index $k$ and this is not possible. In order to avoid a such contradiction, $q^{\ast}E$ is stable, hence $\pr^{\ast}E$ is stable as well.\smallskip

\noindent From all this and by \cite[Proposition 4.5]{IB:GS}, $F$ is a slope stable vector bundle; by construction, $F_{\vert B\times_{\C}\{Q\}}\cong\End(E)$ is stable for $Q\in C_g(\C)$ general.
\end{example}

\subsection{Nilpotent, H-nflat Higgs bundles over Jacobian elliptic surfaces}

\noindent Here we assume that $X$ is a Jacobian elliptic surface, without any other assumption, and we prove Conjecture \ref{conj} for these surfaces (cfr. also Theorem \ref{th3.1}).
\begin{theorem}\label{th4.1}
The Chern classes of any H-nflat Higgs bundle over $X$ are trivial.
\end{theorem}
\noindent This theorem is the second main result of this paper; in order to prove it, we premise some useful properties.
\begin{proposition}\label{prop4.3}
Let $f\colon Y\to Z$ be a flat morphism of schemes, with $Z$ reduced and of finite type over $\spec(\F)$, where $\F$ is a field. Let $g\colon V_1\to V_2$ be a morphism of vector bundles over $Y$. If $g$ restricted to each closed fibre of $f$ is the zero morphism, then $g$ is the zero morphism.
\end{proposition}
\begin{prf}
Let $\spec(B)\subseteq Z$ and $\spec(A)\subseteq f^{-1}(\spec(B))\subseteq Y$ be open affine subschemes, such that $V_1$ and $V_2$ are both trivial vector bundles over $\spec(A)$. Without loss of generality, one may assume that $f_{\vert\spec(A)}$ is surjective, \emph{i.e.} $f^{\sharp}(\spec(A))\colon B\to A$ is a faithfully flat morphism (\cite[\href{https://stacks.math.columbia.edu/tag/00HQ}{tag 00HQ}]{TSP}), and $B$ is a $\F$-algebra of finite type. Let $M_1$ and $M_2$ be free $A$-modules of finite rank such that $V_{k\vert\spec(A)}$ is the sheafification of $M_k$; let $g_0\colon M_1\to M_2$ be a morphism of $A$-modules such that $g_{\vert\spec(A)}$ is the sheafification of $g_0$. By the hypothesis, $g_0\otimes_BB/\fm=0$ for any maximal ideal $\fm$ of $B$; $\displaystyle g_0\left(M_1\right)\subseteq\bigcap_{\fm\in\specm(B)}\fm M_2$, equivalently. By faithfully flatness of $f^{\sharp}$, $M_2$ is a free $B$-module of finite rank (cfr. \cite[Definition 4.2.10, Remark 4.3.1.(ii), Proposition 4.4.1.(i) and (v)]{B:S}). Reasoning via a finite basis of $M_2$ viewed as a $B$-module one has $\displaystyle g_0\left(M_1\right)\subseteq\left(\bigcap_{\fm\in\specm(B)}\fm\right)M_2=\rJ(B)M_2$, where $\rJ(B)$ is the \emph{Jacobson radical of} $B$. By \cite[Corollary 3.2.5]{B:S}, $\rJ(B)=\rad(B)$ (the \emph{nilradical of} $B$) which is $0$, by reducedness of $Z$. Thus $g_0\left(M_1\right)\subseteq\rad(B)M_2=0$ hence $g$ is the zero morphism.~\end{prf}
\begin{remark}
More in general, the previous proposition works if $Z$ is \emph{Jacobson} instead to be of finite type over $\spec(\F)$ for some field $\F$, \emph{i.e.} there exists an affine open covering $\left\{\spec\left(R_i\right)\right\}_{i\in I}$ of $Z$ such that $J\left(R_i\right)=\rad\left(R_i\right)$ for each $R_i$ (cfr. \cite[\href{https://stacks.math.columbia.edu/tag/01P4}{tag 01P4}]{TSP}).
\end{remark}
\begin{lemma}\label{lem4.3}
Let $\fE=(E,\varphi)$ be a rank $r\geq2$, nilpotent, fibrewise semistable Higgs bundle over $X$ with $c_1(E)=0$. Then $\varphi=0$.
\end{lemma}
\noindent One may observe that this statement is blatantly trivial for any nilpotent Higgs line bundle.
\begin{prf}
Let $r=\rank(E)\geq2$ and let $b\in B(\C)$. One has to distinguish two eventualities
\begin{itemize}
\item the fibre $X_b$ is smooth,
\item the fibre $X_b$ is singular.
\end{itemize}
Let $X_b$ be smooth, by \cite[Theorem 4.14]{F:GP:N} $\Gr\left(\fE_{\vert X_b}\right)$ is the direct sum of $r$ stable, H-nflat Higgs line bundles. This is equivalent to the existence of a filtration of $\fE_{\vert X_b}$ in $r$ Higgs subbundles $\fE_s$, each one is H-nflat of rank $s$ (see the proof of \cite[Theorem 3.2]{B:C}). By definition $\varphi_{\vert X_b}\left(\fE_1\right)\subseteq\fE_1\otimes\Omega^1_{X\vert X_b}$ and by nilpotence assumption it turns out that $\varphi_{\vert X_b}\left(\fE_1\right)=0$. Let $r=2$ then $\fE_{\vert X_b}$ fits in a short exact sequence of stable H-flat Higgs line bundles (cfr. \cite[Lemma 3.1]{B:C})
\begin{displaymath}
\xymatrix{
0\ar[r] & \fE_1=\left(\ker\left(\varphi_{\vert X_b}\right),0\right)\ar@{^{(}->}[r] & \fE_{\vert X_b}\ar@{->>}[r] & \fQ=\left(E_{\vert X_b}/\ker\left(\varphi_{\vert X_b}\right),\wtvarphi\right)\ar[r] & 0
};
\end{displaymath}
reasoning on the stalks:
\begin{displaymath}
\forall x\in X_b,s\in E_{\vert X_b,x},\,\varphi_{\vert X_b,x}(s)\in\ker\left(\varphi_{\vert X_b}\right)_x,
\end{displaymath}
this holds by nilpotence of $\varphi_{\vert X_b}$. Thus $\wtvarphi=0$ hence $\varphi_{\vert X_b}=0$. Let $r=3$, repeating the previous reasoning, $\fE_{\vert X_b}$ is an extension of H-nflat Higgs bundles of rank $1$ and $2$, respectively; by the previous step, it follows $\varphi_{\vert X_b}$ vanishes. Assume that all nilpotent, H-nflat Higgs bundle over $\fE_{\vert X_b}$ of rank less or equal than $r$ have Higgs field equals to zero; let $\fE_{\vert X_b}$ have rank $r+1$, then it fits in a short exact sequence of nilpotent, H-nflat Higgs bundles of rank less or equal than $r$. By inductive hypothesis, the Higgs field of $\fE_{\vert X_b}$ vanishes hence, by \emph{Induction Principle}, the Higgs field of nilpotent, H-flat Higgs bundles over $X_b$ vanishes.\smallskip

\noindent Let $X_b$ singular, by Remark \ref{rem4.1} and \cite[table I.4.1]{M:R}, the normalisation $\left(\wtC,\nu\right)$ of $X_b$ is a union of smooth rational curves $\wtC_1,\dotsc,\wtC_m$. By \cite[Th\'eor\`eme 2.1]{G:A} $\Gr\left(\wfE_k\right)=\wfE_k=\left(\nu^{\ast}\fE_{X_b}\right)_{\vert\wtC_k}$ is the rank $r$ trivial bundle over $\wtC_k$ for any $k\in\{1,\dotsc,m\}$. Repeating the previous reasoning, one proves that each $\wfE_k$ is the trivial extension of the trivial Higgs bundle iterated $r$ times. Thus $\varphi_{\vert X_b}$ induces $0\colon\wfE_k\to\wfE_k\otimes_{\cO_{\wtC_k}}\left(\nu^{\ast}\Omega^1_{X\vert X_b}\right)_{\vert\wtC_k}$. Let $\wteta_k$ be the generic point of $\wtC_k$ and let $\eta_1,\dotsc,\eta_1$ be the generic points of $X_b$. By definition
\begin{displaymath}
\wfE_{\wteta_k}=\fE_{\eta_k}\otimes_{\cO_{X_b,\eta_k}}\cO_{\wtC_k,\wteta_k},\,\left(\nu^{\ast}\Omega^1_{X_b}\right)_{\wteta_k}=\left(\Omega^1_{X\vert X_b}\right)_{\eta_k}\otimes_{\cO_{X_b,\eta_k}}\cO_{\wtC_k,\wteta_k},
\end{displaymath}
for any $k\in\{1,\dotsc,m\}$. Hence, up to canonical isomorphism
\begin{displaymath}
\left(\varphi_{\vert X_b}\otimes\Id\right)_{\wteta_k}=0_{\wteta_k}\colon\fE_{\eta_k}\otimes_{\cO_{X_b,\eta_k}}\cO_{\wtC_k,\wteta_k}\to\left(\fE_{\eta_k}\otimes_{\cO_{X_b,\eta_k}}\left(\Omega^1_{X\vert X_b}\right)_{\eta_k}\right)\otimes_{\cO_{X_b,\eta_k}}\cO_{\wtC_k,\wteta_k};
\end{displaymath}
since $\cO_{\wtC_k,\wteta_k}=\kappa\left(\wteta_k\right)$ is a finite degree field extension of $\cO_{X_b,\eta_k}=\kappa\left(\eta_k\right)$, $\_\otimes_{\kappa\left(\eta_k\right)}\kappa\left(\wteta_k\right)$ is a faithfully flat morphism, hence the vanishing of $\left(\varphi_{\vert X_b}\otimes\Id\right)_{\wteta_k}$ implies the vanishing of $\varphi_{\eta_k}$. Equivalently, any $\Im\left(\varphi_{\vert\overline{\left\{\eta_k\right\}}}\right)$ is a torsion subsheaf of $\Omega^1_{X\vert\overline{\left\{\eta_k\right\}}}$; since this last one is a locally free sheaf then $\Im\left(\varphi_{\vert\overline{\left\{\eta_k\right\}}}\right)=\uzero_{\vert\overline{\left\{\eta_k\right\}}}$, \emph{i.e.} $\varphi_{\vert\overline{\left\{\eta_k\right\}}}=0$. Since all this holds on any irreducible component of $X_b$, it turns out that $\varphi_{\vert X_b}$ vanishes.\smallskip

\noindent By the previous proposition, one has the claim.
\end{prf}
\noindent \textbf{Proof of Theorem \ref{th4.1}.}
Let $\fE=(E,\varphi)$ be a rank $r$ H-nflat Higgs bundle over $X$. If $r=2$ then $\cS^{p,\cL}_X$ is well-defined by Lemma \ref{lem3.2} (see also Proposition \ref{prop3.2}). Then $\displaystyle\lim_{z\to0}[\fE]_{\sim}=\left[\fE_0\right]_{\sim}$ is the JH-equivalent class of a nilpotent, semistable, H-nflat Higgs bundle $\fE_0=\left(E_0,\varphi_0\right)$. By the previous lemma, $\varphi_0=0$ and therefore
\begin{displaymath}
P(t)=\int_X\frac{1}{6}\left[c_1(X)^2+c_2(X)\right]+c_1(L)c_1(X)t+c_1(L)^2t^2\in\Q[t]
\end{displaymath}
by \cite[Proposition 1.3]{C:P}, \emph{i.e.} the Chern classes of $\fE$ are trivial. Thus one repeats the same reasoning for rank $3$ H-nflat Higgs bundles and proves the triviality of Chern classes of these Higgs bundles. By induction, one may assume that the rank $r$ H-nflat Higgs bundles have trivial Chern classes, hence one repeats again the same reasoning and proves the triviality of rank $r+1$ H-nflat Higgs bundles. By \emph{Induction Principle}, one has the claim.\hfill\text{(Q.e.d.)}\medskip

\noindent By Lemma \ref{lem3.2} it follows the following corollary.
\begin{corollary}
H-nflat Higgs bundles over $X$ are semistable.
\end{corollary}
\begin{corollary}[{cfr. \cite[Corollary 2.11]{B:LG}}]
A Higgs bundle $\fE=(E,\varphi)$ over $X$ is H-nflat if and only if $E$ is nflat.
\end{corollary}
\begin{prf}
By Remark \ref{rem3.1}.\ref{rem3.1.b} and the previous corollary, one needs to prove the ``only if'' part. Let $\fE$ be H-nflat and let $[\fE]_{\sim}$ be its JH-equivalence class in the moduli space $\rM^{p,\cL}_X$; where $p$ is an opportune Hilbert polynomial with respect to $\cL$, an ample invertible sheaf on $X$. Since $\fE$ is curve semistable, by Lemma \ref{lem4.3} and Remark \ref{remC.1}, $\displaystyle\lim_{z\to0}z\cdot[\fE]_{\sim}=[(E,0)]_{\sim}$, by Corollary \ref{cor3.1} $E$ is nflat.
\end{prf}
\begin{corollary}
A Higgs bundle $\fE=(E,\varphi)$ over $X$ is curve semistable if and only if $E$ is a curve semistable vector bundle.
\end{corollary}
\begin{prf}
If $\fE$ is curve semistable then $\End(\fE)$ is H-nflat (cfr. Proposition \ref{prop3.1}). By the previous corollary, $\End(E)$ is nflat, equivalently; reversing this reasoning one has the claim.
\end{prf}

\appendix

\section{Moduli stack of stable maps}\label{ModStackStMaps}

\subsection{Families of $m$-pointed prestable curves}

\noindent We recall that proper algebraic spaces of dimension at most $1$ over a field are projective scheme (\cite[\href{https://stacks.math.columbia.edu/tag/0A26}{tags 0A26} and \href{https://stacks.math.columbia.edu/tag/0AE7}{0AE7}]{TSP})
. This property will be used in the next, tacitly.
\begin{definition}
A connected, nodal, projective curve $C$ of arithmetic genus $g$ with an ordered collection $\left(p_1,\dotsc,p_m\right)$ of $m$, pairwise distinct, closed, smooth points is \emph{prestable} if $(g,n)\neq(1,0)$.\smallskip

\noindent Nodal and marked points of $\left(C,p_1,\dotsc,p_m\right)$ are called \emph{special points}.
\end{definition}
\begin{definition}
By \emph{family of} $m$-\emph{pointed prestable curves} $\cC$ \emph{over a scheme} $S$ we mean a flat, proper and finitely presented morphism $\cC\to S$ of \emph{algebraic spaces} together with $m$ sections $\sigma_1,\dotsc,\sigma_m\colon S\to\cC$, such that every fibre $\left(\cC_s,\sigma_1(s),\dotsc,\sigma_m(s)\right)$ is a prestable curve of genus $g$ over $\spec(\kappa(s))$.
\end{definition}
\noindent One defines families of $m$-pointed curves of genus $g$ analogously. Where it is possible, we shall skip any reference to the genus. We shall write $\left(\cC\to S,\sigma_{\ast}\right)$ instead of $\left(\cC\to S,\sigma_1,\dotsc,\sigma_m\right)$, for short.
\begin{proposition}[cfr. {\cite[Proposition 6.4.1]{A:J}}]\label{propA.1}
Let $\cC\to S$ be a family of $m$-pointed prestable curves. There exists an \'etale cover $\left\{\wtS_t\to S\right\}_{t\in T}$ such that $\wtcC_t=\cC\times_S\wtS_t\to\wtS_t$ is a projective family of $m$-pointed prestable curves for each $t\in T$.
\end{proposition}
\begin{prf}
For a while, let assume $\cC$ be a scheme and let assume $S$ be an affine scheme. By hypothesis, $S$ is a quasi-compact and quasi-separated scheme over $\spec(\Z)$. By \emph{Absolute Noetherian Approximation} (\cite[\href{https://stacks.math.columbia.edu/tag/0GS1}{tags 0GS1} and \href{https://stacks.math.columbia.edu/tag/09MV}{09MV}]{TSP}), there exists a direct set of indexes $I$ and an inverse system $\left(S_i,f_{ij}\right)_{i,j\in I}$ of finitely presented schemes $S_i$ with affine transition morphisms $f_{ij}\colon S_i\to S_j$ over $\spec(\Z)$, such that $\displaystyle\lim_{\overleftarrow{i\in I}}S_i=S$. By \emph{Descent of Morphisms under Limits} (\cite[\href{https://stacks.math.columbia.edu/tag/01ZM}{tag 01ZM}]{TSP}), there exist an index $i_0\in I$ and a finitely presented morphism $\cC_{i_0}\to S_{i_0}$ of schemes such that $\cC\cong\cC_{i_0}\times_{S_{i_0}}S$; moreover, if one sets $\cC_{\ell}=\cC_{i_0}\times_{S_{i_0}}S_{\ell}$ with $\ell\geq i_0\in I$, then $\displaystyle\cC=\lim_{\overleftarrow{\ell\geq i_0\in I}}\cC_{\ell}$. Thus, for $\ell\gg i_0\in I$, $\cC_{\ell}\to S_{\ell}$ is a family of $m$-pointed curves (\cite[\href{https://stacks.math.columbia.edu/tag/04AI}{tags 04AI}, \href{https://stacks.math.columbia.edu/tag/05M5}{05M5}]{TSP}, \cite[Th\'eor\`eme 8.10.5.(vi) and (xii)]{EGA:IV}). Indeed, the curves $\cC_s$ are \emph{geometrically connected} for any geometric point $s\in S$ (\cite[\href{https://stacks.math.columbia.edu/tag/0387}{tag 0387}]{TSP}), hence each of their base changes is connected as well. Let $\sigma\colon S\to\cC$ be a section of this family, one has the following diagram
\begin{displaymath}
\xymatrix{
S_{\ell}\ar@/_1pc/[ddr]_{\Id_{S_{\ell}}}\ar@/^1pc/[drr]^{\sigma\circ f}\\
 & \cC_{\ell}\ar[d]\ar[r] & \cC\ar[d]\\
 & S_{\ell}\ar[r]_f & S
},
\end{displaymath}
by the Universal Property of fibre products of schemes, there exists a unique morphism ${\sigma_{\ell}\colon S_{\ell}\to\cC_{\ell}}$ which makes commutative the whole of diagram, \emph{i.e.} $\sigma_{\ell}$ is a section of the family of curves $\cC_{\ell}\to S_{\ell}$, by \cite[\href{https://stacks.math.columbia.edu/tag/0C5B}{tag 0C5B}]{TSP} this is a family of curves at worst nodal. If the curves parameterised by $S$ have arithmetic genus $1$, then by \cite[\href{https://stacks.math.columbia.edu/tag/0BY9}{tag 0BY9}]{TSP} the curves parameterised by $S_{\ell}$ have geometric genus $1$; and by the previous reasoning, $\cC_{\ell}$ is family of $m$-pointed prestable curves.\smallskip

\noindent Since \'etale and projective properties of morphisms of schemes are stable under base changes, one may replace $\cC\to S$ with $\cC_{\ell}\to S_{\ell}$, where $\ell\gg i_0\in I$. Let $s\in S$, define $S_{\ell}=\spec\left(\cO_{S,s}/\fm^{\ell+1}\right)$ for $\ell\in\N_{\geq0}$ and $\whS=\spec\left(\whcO_{S,s}\right)$. Consider the Cartesian diagram
\begin{displaymath}
\xymatrix{
\cC_s=\cC_0\ar[d]\ar@{^{(}->}[r] & \cC_1\ar[d]\ar@{^{(}->}[r] & \dotsc\ar@{^{(}->}[r] & \whC\ar[d]\ar@{^{(}->}[r] & \cC\ar[d]\\
s=\spec(\kappa(s))\ar@{^{(}->}[r] & S_1\ar@{^{(}->}[r] & \dotsc\ar@{^{(}->}[r] & \whS\ar@{^{(}->}[r] & S
}.
\end{displaymath}
Since $\cC_0$ is a proper, one-dimensional scheme over a field, it is a projective; let $L_0$ be an ample line bundle over $\cC_0$, by \cite[Proposition C.2.11.(3)]{A:J}, the space of the obstruction to deforming $L_0$ to a line bundle $L_{n+1}$ over $C_{n+1}$ is $\rH^2\left(\cC_0,\cO_{\cC_0}\otimes_{\kappa(s)}\fm^n\right)$ which is $0$ since $\dim\cC_0=1$. Thus, there exists a sequence of compatible line bundles $L_n$ over $\cC_n$; by \emph{Grothendieck's Existence Theorem} (\cite[Theorem C.5.3]{A:J}), there exists a line bundle $\whL$ over $\whC$ extending all $L_n$'s. Applying \emph{Artin Approximation Theorem} (\cite[Theorem B.5.18]{A:J}) to the contravariant functor
\begin{gather*}
\left(\sfSch_{/S}\right)^{op}\to\sfSet\\
(T\to S)\mapsto\Pic\left(\cC_T\right)
\end{gather*}
one has an \'etale morphism $\left(\wtS_0,\wts\right)\to(S,s)$ and a line bundle $\wtL_0$ over $\wtcC_0=\cC\times_S\wtS_0$ extending $L_0$. By \cite[Theorem 1.2.17]{L:RK}, there exists an open neighbourhood $\wtS$ of $\wts$ such that $\wtL_{\vert\cC\times_S\wtS}$ is fibrewise ample hence it is ample (\cite[Theorem 1.7.8]{L:RK}). One gets $\wtcC=\cC\times_S\wtS$ and has the claim under the previous assumptions on $\cC$ and $S$.\smallskip

\noindent Let us return to assuming that $S$ is a scheme. Let $s\in S$, let $U$ be an affine open neighbourhood of $s$ and let $i\colon U\hookrightarrow S$ be the open inclusion. By the previous part, there exists an \'etale morphism $\left(\wtS,\wts\right)\to(U,s)\xrightarrow{i}(S,s)$ which satisfies the claim.\smallskip

\noindent Let us return to assuming that $\cC$ is an algebraic space. By \cite[\href{https://stacks.math.columbia.edu/tag/09YC}{tag 09YC}]{TSP}, there exists a commutative diagram
\begin{displaymath}
\xymatrix{
\ocC\ar[dr]\ar[r] & \cC\ar[d]\\
& S
}
\end{displaymath}
where $\ocC$ is a finite $S$-scheme finitely presented over $\cC$ which surjects onto $\cC$. In particular $\ocC$ is proper over $S$, and its fibres are curves. Repeating the previous reasoning, one has the following commutative diagram
\begin{displaymath}
\xymatrix{
\ocC\times_S\wtS\ar[dr]\ar[r] & \wtcC\ar[d]\\
& \wtS
}
\end{displaymath}
where $\ocC\times_S\wtS$ is a projective scheme, and $\wtcC$ is a family of $m$-pointed prestable curves. Since $\ocC\times_S\wtS$ is a separated scheme, such that every finite subset is contained in an open affine subscheme; by \cite[Theorem 3.3]{G:P} this is equivalent to state that every finite set of points of $\left|\wtcC\right|$ is contained in a Zariski-open neighbourhood that is representable by an affine scheme, \emph{i.e.} $\wtcC$ is a scheme and the proof is completed.
\end{prf}

\subsection{Families of stable maps}

\noindent Let $X$ be a smooth, projective $\K$-variety.

\begin{definition}\label{defA.1}
Let $\beta\in\rA^{n-1}(X)$. A \emph{stable map into} $X$ is the data $\left(C,p_1,\dotsc,p_m,f\right)$ of a $m$-pointed prestable curve with a morphism $f\colon C\to X$ such that $f_{\ast}[C]=\beta$ and every irreducible component $C_0$ of $C$ of arithmetic genus $0$ ($1$, respectively) that is contracted to a point under $f$, $C_0$ has at least $3$ ($1$, respectively) special points.
\end{definition}
\begin{definition}
A \emph{family of stable maps into} $X$ \emph{over a scheme} $S$ is a tuple $\left(\cC\to S,\cC\to X,\sigma_{\ast}\right)$ such that $\left(\cC\to S,\sigma_{\ast}\right)$ is a family of $m$-pointed prestable curve over $S$, and $\cC\to X$ restricted to geometric fibres of $\cC\to S$ is a stable map.
\end{definition}
\begin{remark}\label{remA.1}
To give a diagram
\begin{displaymath}
\xymatrix{
\cC\ar[d]\ar[r] & X\\
S\ar@/^1pc/[u]^{\sigma_i}
}
\end{displaymath}
is equivalent to give a morphism $\cC\to X\times_{\Z}S$ with $m$ sections $\sigma_i\circ\pr_S$.
\end{remark}
\noindent Let $\osfM_{g,m}(X,\beta)$ be the category of families of stable maps into $X$ with fixed class $\beta\in\rA^{n-1}(X)$; here a morphism is the following commutaive diagram
\begin{displaymath}
\xymatrix{
\cC_1\ar[d]\ar[r]_g\ar@/^1pc/[rr] & \cC_2\ar[d]\ar[r] & X\\
S_1\ar[r]_f\ar@/^1pc/[u]^{\sigma_{i,1}} & S_2\ar@/_1pc/[u]_{\sigma_{i,2}}
},
\end{displaymath}
\noindent where:
\begin{enumerate}[a)]
\item $\cC_1\to X$ and $\cC_2\to X$ restricted to geometric fibres are stable maps with fixed class $\beta$;
\item the square is Cartesian in the category of algebraic spaces $\sfAlgSp$,
\item $\sigma_{2,i}\circ f=g\circ\sigma_{1,i}$ for any $i\in\{1,\dotsc,m\}$.
\end{enumerate}
\noindent The composition of morphisms is defined juxtaposing these diagrams in the obvious way. 
\begin{remark}\label{remA.2}
In the previous diagram, if one forgets $X$ then one has the definition of morphisms of families of $m$-pointed prestable curves over schemes. 
\end{remark}
\begin{proposition}[cfr. {\cite[Proposition 6.4.1]{A:J}}]\label{propA.2}
Let $\left(\cC\to S,\cC\to X,\sigma_{\ast}\right)$ be a family of stable maps into $X$. There exists an \'etale cover $\left\{\wtS_t\to S\right\}_{t\in T}$ such that $\left(\wtcC_t=\cC\times_S\wtS_t\to\wtS_t,\wtcC_t\to X,\wtsigma_{\ast}\right)$ is a projective family of stable maps into $X$ for each $t\in T$.
\end{proposition}
\begin{prf}
By Proposition \ref{propA.1} there exists a such \'etale cover for families of $m$-pointed stable curves. Let $s_0\in S$ be a geometric point, let $C_0$ be an irreducible components of $\cC_{s_0}$ of arithmetic genus $0$ ($1$, respectively). If $C_0$ contracted to a pointed then it has at least $3$ ($1$, respectively) special points. Let $\wts_0\in\wtS$ be a geometric point over $s_0$; by \cite[\href{https://stacks.math.columbia.edu/tag/0BY9}{tags 0BY9} and \href{https://stacks.math.columbia.edu/tag/0C57}{0C57}]{TSP}, the base change $\wtC_0$ of $C_0$ over $\wts_0$ is a projective curve, of geometric genus $0$ ($1$, respectively) and it is at worst nodal; if $C_0$ has also marked points, by construction $\wtsigma_{\ast}$'s determine marked points on $\wtC_0$. By all this, $\wtC_0$ has the same number of special points of $C_0$ and one has the claim.
\end{prf}
\begin{theorem}[{\cite[Theorem in subsubsection 1.3.1]{K:M}}]\label{thA.1}
$\osfM_{g,m}(X,\beta)$ is a proper, algebraic stack with finite diagonal.
\end{theorem}

\section{Projective varieties and base change}

\noindent We fix the following definitions.
\begin{definitions}\label{defB.1}
Let $\F$ be a field and let $\F^{\prime}$ be an its extension field.
\begin{enumerate}[a)]
\item One says $\F^{\prime}$ \emph{algebraic separable} if it is algebraic and the minimal polynomial $p$ of any element of $\F^{\prime}$ is \emph{coprime} with its derivative $Dp$ (\cite[\href{https://stacks.math.columbia.edu/tag/09H1}{tag 09H1}]{TSP}).
\item\label{defB.1.b} A collection of elements $\{x_i\}_{i\in I}$ of $\F^{\prime}$ is called \emph{algebraically independent} over $\F$ if the map
\begin{gather*}
\F[X_i\mid i\in I]\to\F^{\prime}\\
X_i\mapsto x_i
\end{gather*}
is injective (\cite[\href{https://stacks.math.columbia.edu/tag/030E}{tag 030E}]{TSP}).
\item A \emph{transcendence basis of $\F^{\prime}$ over $\F$} is a collection of elements $\{x_i\}_{i\in I}$ which are algebraically independent over $\F$ and such that $\F^{\prime}$ is an algebraic field extension of $\F(x_i\mid i\in I)$ (\cite[\href{https://stacks.math.columbia.edu/tag/030E}{tag 030E}]{TSP}).
\item One says $\F^{\prime}$ is \emph{separably generated} over $\F$ if there exists a transcendence basis $\left\{x_i\in\F^{\prime}\right\}_{i\in I}$ such that $\F^{\prime}$ is an algebraic separable field extension of $\F(x_i\mid i\in I)$ (\cite[\href{https://stacks.math.columbia.edu/tag/030O}{tag 030O}]{TSP}).
\item One says $\F^{\prime}$ \emph{separable} over $\F$ if for any field extension $\F\subseteq\F^{\second}\subseteq\F^{\prime}$, with $\F^{\second}$ finitely generated over $\F$, $\F^{\second}$ is separably generated over $\F$ (\cite[\href{https://stacks.math.columbia.edu/tag/030O}{tag 030O}]{TSP}).
\item One says $\F$ \emph{perfect} if any its field extension is separable (\cite[\href{https://stacks.math.columbia.edu/tag/03=Y}{tag 030Y}]{TSP}).
\end{enumerate}
\end{definitions}
\begin{remarks}\label{remB.2}
\,
\begin{enumerate}[a)]
\item Any field extension has a transcendence basis (\cite[\href{https://stacks.math.columbia.edu/tag/030F}{tag 030F}]{TSP}).
\item Any separably generated field extension is separable (\cite[\href{https://stacks.math.columbia.edu/tag/030X}{tag 030X}]{TSP}).
\item\label{remB.2.c} The algebraic field extensions of a characteristic $0$ field are separable, hence any field of characteristic $0$ is perfect. 
\item A field of characteristic $0$ is separably closed if and only if it is algebraically closed. \hfill{$\Diamond$}
\end{enumerate}
\end{remarks}
\noindent From now on, let $X$ be a smooth, projective variety defined over $\K$, let $\F$ be a field extension of $\K$ and let $X_{\F}=X\times_{\K}\spec(\F)$. By \cite[\href{https://stacks.math.columbia.edu/tag/01WF}{tags 01WF} and \href{https://stacks.math.columbia.edu/tag/020J}{020J}]{TSP} and Remark \ref{remB.2}.\ref{remB.2.c}, $X_{\F}$ is an irreducible, smooth, projective scheme of finite type over $\F$. Moreover, one has the following lemma.
\begin{lemma}\label{lemB.1}
$X_{\F}$ is a smooth projective variety of dimension $\dim(X)$.
\end{lemma}
\begin{prf}
$X$ is a reduced scheme \emph{i.e.} for any open affine subset $U$ of $X$, $\cO_X(U)$ is a reduced $\K$-algebra hence it is a \emph{geometrically reduced $\K$-algebra} (\cite[\href{https://stacks.math.columbia.edu/tag/030S}{tags 030S} and \href{https://stacks.math.columbia.edu/tag/030U}{030U}]{TSP}). Since $\K$ is a perfect field, then $\F$ is a separable field extension of $\K$ hence $\cO_X(U)\otimes_{\K}\F$ is a reduced $\F$-algebra (\cite[Proposition 10.5.22]{FOAG}). Consider the following Cartesian diagram
\begin{displaymath}
\xymatrix{
X_{\F}\ar[r]^{f}\ar[d] & X\ar[d]\\
\spec(\F)\ar[r] & \spec(\K)
}.
\end{displaymath}
From all this $\left\{f^{-1}(U)=\spec\left(\cO_X(U)\otimes_{\K}\F\right)\right\}_{\substack{U\subseteq X\\
\text{open and affine}}}$ is an affine, open covering of $X_{\F}$ given by reduced subschemes,\emph{i.e.} $X_{\F}$ is a reduced scheme and one infers that $X_{\F}$ is an integral scheme. The statement on the dimension of $X_{\F}$ follows from \cite[Exercise 12.2.M]{FOAG}.
\end{prf}
\begin{remark}\label{remB.1}
By the previous proof, $f$ is a \emph{quasi-affine morphism}, \emph{i.e.} there exists an open affine covering $\left\{U_i\right\}_{i\in I}$ of $X$ such that $f^{-1}\left(U_i\right)$ is an open subscheme of an affine scheme for any $i\in I$.~\end{remark}
\begin{definition}[{see \cite[Subsection A.2.3]{A:J}}]\label{defB.2}
A \emph{fpqc}\footnote{
In French, ``\emph{fid\`element plat et quasi-compact}''.
} morphism $\alpha\colon Y\to Z$ of schemes is a faithfully flat morphism\footnote{
In other words, $\alpha$ is a surjective flat morphism of schemes (see \cite[Subsection A.2.3]{A:J}).
} for which there exists an affine open covering $\{Z_i\}_{i\in I}$ of $Z$, such that any $Z_i$ is the image of a quasi-compact open subset of $Y$.
\end{definition}
\begin{proposition}\label{propB.1}
The canonical morphism $f\colon X_{\F}\to X$ is fpqc.
\end{proposition}
\begin{prf}
Since $X_{\F}$ is a closed subscheme of $\bP^N_{\F}$ (\cite[Proposition 7.3.13]{B:S}) for some $N\in\N_{\geq1}$, $X_{\F}$ is quasi-compact. Thus one may consider a finite, affine, open covering $\left\{f^{-1}(U_i)\right\}_{i\in\{1,\dotsc,m\}}$ of $X_{\F}$, where each $U_i$ is an open, affine subscheme of $X$. By \cite[Proposition 4.4.1.iii and Corollary 7.2.7]{B:S}, $f$ is faithfully flat, $f^{-1}(U_i)=\spec\left(\cO_X(U_i)\otimes_{\K}\F\right)$ for any $i\in\{1,\dotsc,m\}$ and these are quasi-compact topological spaces. In other words ${U_i=f\left(\spec\left(\cO_X(U_i)\otimes_{\K}\F\right)\right)}$ for any $i\in\{1,\dotsc,m\}$, \emph{i.e.} the claim holds.
\end{prf}

\section{Moduli spaces of Higgs sheaves}\label{HiggsModSp}

\noindent Here $X$ be a smooth, complex, projective variety. Let $\fE=(\cE,\varphi)$ be a semistable Higgs sheaf on $X$. There exists a filtration
\begin{displaymath}
\uzero=\fE_0\subsetneqq\fE_1\subsetneqq\dotsc\subsetneqq\fE_{m-1}\subsetneqq\fE_m=\fE
\end{displaymath}
of \emph{saturated}\footnote{
Let $\cF$ be a subsheaf of $\cE$, it is \emph{saturated} if $\cE/\cF$ is a torsion-free sheaf.
} Higgs subsheaves of $\fE$ such that each successive quotient $\fE_k/\fE_{k-1}$ is stable and with Hilbert polynomial $p_{\fE_k/\fE_{k-1},norm}=p_{\fE,norm}$ with respect to an ample invertible sheaf $\cL$ on $X$. This is called \emph{Jordan-H\"older filtration of} $\fE$ (\emph{JH-filtration}, for short); it is not unique, but the \emph{graded Higgs sheaf} $\displaystyle\Gr(\fE)=\bigoplus_{k=1}^m\fE_k/\fE_{k-1}$ is unique up to isomorphisms (\cite[Section 4]{HC:MG}). One declares \emph{JH-equivalent} two semistable Higgs sheaves $\fE_1$ and $\fE_2$ if ${\Gr\left(\fE_1\right)\cong\Gr\left(\fE_2\right)}$.\medskip

\noindent By a \emph{family of semistable Higgs sheaves on} $X$ \emph{parameterised by a} $\C$-\emph{scheme} $S$ we mean a Higgs sheaf $\fF=(\cF,\Psi)$ on $X_S=X\times_{\C}S$ such that $\fF_{\vert X_{\kappa(s)}}$ is a semistable Higgs sheaf and $\cF_{\vert X_{\kappa(s)}}$ is a flat $\cO_{S,s}$-module for any $s\in S$; $X_{\kappa(s)}=X\times_{\C}\spec(\kappa(s))$ of course. Two such families $\fF_1$ and $\fF_2$ are \emph{equivalent} if they are isomorphic as Higgs sheaves.
\begin{remark}
By the assumptions, $X$ is a projective scheme of finite presentation over $\spec(\C)$. Since these properties are stable under base change, one may apply \cite[Proposition 23.150]{G:W:2} which proves the constancy of Hilbert polynomial of the fibres of a family of Higgs sheaf on $X$ parameterised by a scheme $S$ on $\spec(\C)$.
\end{remark}
\noindent Consider the contravariant functor ${\bf M}^{p,\cL}_X\colon\left(\sfSch_{/\C}\right)^{opp}\to\sfSet$ under which
\begin{displaymath}
S\mapsto\left\{\begin{tabular}{l}
families of semistable Higgs sheaves on $X$ parameterised\\ by $S$ whose fibres have $p$ as Hilbert polynomial
\end{tabular}\right\}/\text{equivalence}.
\end{displaymath}
\noindent where $p\in\Q[\lambda]$ and $\deg(p)=\dim(X)$. Simpson has proved that this functor is corepresentable by a quasi-projective scheme $\rM^{p,\cL}_X$ (\cite[Theorem 4.7]{S:CT:2} and \cite[Lemma 6.5]{S:CT:3}). Moreover, there is a bijection between closed points of $\rM^{p,\cL}_X$ and JH-equivalence classes of Higgs sheaves on $X$ with Hilbert polynomial $p$. $\rM^{p,\cL}_X$ is called \emph{moduli space of Higgs sheaves on} $X$ \emph{with fixed Hilbert polynomial} $p$ \emph{with respect to the ample invertible sheaf} $\cL$.
\begin{remark}\label{remC.1}
There exists a holomorphic action of $\C^{\times}$ on $\rM^{p,\cL}_X$ given by
\begin{displaymath}
\cdot\colon\left(z,[(E,\varphi)]_{\sim}\right)\in\C^{\times}\times\rM^{p,\cL}_X\to[(E,z\varphi)]_{\sim}\in\rM^{p,\cL}_X
\end{displaymath}
(see \cite[Subsection 6.5]{S:CT:3}). In particular, the limit $\displaystyle\lim_{z\to0}z\cdot[\fE]_{\sim}$ there exists in $\rM^{p,\cL}_X$ (see proof of \cite[Corollary 6.12]{S:CT:3}). Since these Higgs bundles ``limits'' are fixed points of this action, by \cite[Lemma 4.1]{S:CT:1} they are semistable \emph{nilpotent} Higgs bundles, \emph{i.e.} let $\fE_0=\left(E,\varphi_0\right)$ be a such limit, there exists $E_1,\dotsc,E_m$ vector subbundles of $E$ such that $\displaystyle E=\bigoplus_{k=1}^mE_k$, $\varphi_0\colon E_k\to E_{k+1}\otimes\Omega^1_X$ for any $k\in\{1,\dotsc,m-1\}$ and $\varphi_0\left(E_m\right)=0$.
\end{remark}
\noindent About these ``limits'', We premise the following general proposition.
\begin{proposition}\label{propC.1}
Let $Y$ be a smooth, projective surface defined over $\K$, let $H$ be an ample line bundle over $Y$, let $\cV_1$ and $\cV_2$ locally free sheaves on $Y$ of finite rank, and let $f\colon\cV_1\to\cV_2$ be a morphism of sheaves. If for $m\gg1$ there exists a smooth projective curve $C\in\lvert mH\rvert$ such that $f_{|C}=0$ then $f=0$.
\end{proposition}
\begin{prf}
By \emph{Bertini's Theorem} (\cite[Theorem II.8.18 and Remark III.7.9.1]{H:RC:1}), such $m$ and $C$ exist. Let
\begin{displaymath}
\cHom\left(\cV_1,\cV_2\right)\cong\cV_1^{\vee}\otimes_{\cO_Y}\cV_2=\cW.
\end{displaymath}
Since $C$ is a closed subscheme of $Y$, one has the following short exact sequence
\begin{displaymath}
\xymatrix{
0 \ar[r] & \cO_Y(-mH)\ar[r] & \cO_Y\ar[r] & i_{\ast}\cO_C\ar[r] & 0
}
\end{displaymath}
where $i\colon C\hookrightarrow Y$ is the inclusion. Tensorising it with $\cW$, one has
\begin{displaymath}
\xymatrix{
0\ar[r] & \cW(-mH)\equiv\cW\otimes_{\cO_Y}\cO_Y(-mH)\ar[r] & \cW\ar[r] & i_{\ast}\cW_{\vert C}\cong\cW\otimes_{O_Y}i_{\ast}\cO_C\ar[r] & 0
}.
\end{displaymath}
Passing to long short exact sequence in cohomology, one has
\begin{displaymath}
\xymatrix{
0\ar[r] & H^0\left(Y,\cW(-mH)\right)\ar[r] & H^0(Y,\cW)\ar[r] & H^0\left(C,\cW_{\vert C}\right)\ar[r] & \dotsc
}
\end{displaymath}
By \emph{Serre Duality}
\begin{displaymath}
H^0\left(Y,\cW(-mH)\right)=H^2\left(Y,\cW^{\vee}(mH)\otimes_{\cO_Y}\omega_Y\right)^{\vee},
\end{displaymath}
by \emph{Serre Vanishing Theorem} the right hand side vanishes hence the claim holds.
\end{prf}
\begin{corollary}\label{corC.1}
Let $\fE=(E,\varphi)$ be a Higgs bundle over $X$ such that $E$ is slope semistable with respect to the polarization determined by $\cL$. Then
\begin{displaymath}
\lim_{z\to0}z\cdot[\fE]_{\sim}=[(E,0)]_{\sim},
\end{displaymath}
where the notations are obvious.
\end{corollary}
\begin{prf}
Let $H$ be the ample Cartier divisor of $X$ such that $\cO_X(H)=\cL$. By \emph{Bertini's Theorem} (\cite[Theorem II.8.18 and Remark III.7.9.1]{H:RC:1}), there exist $m\gg1$ and $C\in\lvert mH\rvert$ such that $C$ is a smooth, projective curve. By \cite[Theorem 6.1]{MVB:RA}, $E_{\vert C}$ is semistable; repeating the proof of Proposition \ref{prop3.2}, one has:
\begin{displaymath}
\left[\left(E_{\vert C},\left(\varphi_0\right)_{\vert C}\right)\right]_{\sim}=i^{\sharp}\left(\left[(E,\varphi_0)\right]_{\sim}\right)=i^{\sharp}\left(\lim_{z\to0}z\cdot[(E,\varphi)]_{\sim}\right)=\lim_{z\to0}z\cdot\left[\left(E_{\vert C},\varphi_{\vert C}\right)\right]_{\sim}=\left[\left(E_{\vert C},\left(\varphi_{\vert C}\right)_0\right)\right]_{\sim}
\end{displaymath}
where $i\colon C\hookrightarrow X$ is the inclusion. By \cite[Proposition 2.2]{G:Z-R} and \cite[Theorem C]{F:Y}, ${\left(\varphi_0\right)_{\vert C}=\left(\varphi_{\vert C}\right)_0=0}$; by the previous proposition, $\varphi_0=0$.
\end{prf}
\bigskip

\noindent{\bf Acknowledgment.} I would like to thank in special way Jason Starr and the anonymous referee who provided insightful suggestions on the proofs of Lemma \ref{lem2.2} and Corollary \ref{cor2.1}. Proposition \ref{prop4.1} was brought to my attention by S\'andor Kova\'cs and Angelo Vistoli. The idea of the proof of Proposition \ref{prop4.3} is due to Alexander Kuznetsov. Proposition \ref{propC.1} and Corollary \ref{corC.1} have been achieved after some discussion with Damian R\"ossler and Ronald Alberto Z\'u\~niga-Rojas, respectively. I thank all of them warmly. I am very grateful to my Ph.D. advisors Ugo Bruzzo and Beatriz Gra\~na Otero for their help, their energy and their support. I thank S.I.S.S.A. for the hospitality while part of this work was done.\medskip

\noindent{\bf Statement about competing or financial interests.} The author has no competing or financial interests to declare that are relevant to the content of this article.

\end{document}